\documentclass{elsarticle}
\usepackage{amssymb,bm}
\usepackage{subcaption}
\usepackage{hyperref}
\usepackage{amsmath}
\usepackage{algorithm2e}
\usepackage{graphicx}

\usepackage{tikz}
\newcommand*\circled[1]{\tikz[baseline=(char.base)]{
            \node[shape=circle,draw,inner sep=1pt] (char) {#1};}}

\journal{Journal of Computational Physics}
\begin{document}

\begin{frontmatter}

\title{A High Resolution PDE Approach to Quadrilateral Mesh Generation}

\author[add1]{Julian Marcon}
\ead{julian.marcon14@imperial.ac.uk}
\author[add2]{David A. Kopriva\corref{cor1}}
\ead{kopriva@math.fsu.edu}
\author[add1]{Spencer J. Sherwin}
\ead{s.sherwin@imperial.ac.uk}
\author[add1]{Joaquim Peir\'{o}}
\ead{j.peiro@imperial.ac.uk}

\address[add1]{
  Imperial College London,
  South Kensington Campus,
  London SW7 2AZ, UK
}
\address[add2]{
  The Florida State University,
  Tallahassee, FL 32306, USA
  and
  San Diego State University,
  San Diego, CA 92182, USA
}

\cortext[cor1]{Corresponding author}

\begin{abstract}

We describe a high order technique to generate quadrilateral decompositions and meshes for complex two dimensional domains using spectral elements in a field guided procedure.
Inspired by cross field methods, we never actually compute crosses.
Instead, we compute a high order accurate guiding field using a continuous Galerkin (CG) or discontinuous Galerkin (DG) spectral element method to solve a Laplace equation for each of the field variables using the open source code \emph{Nektar++}.
The spectral method provides spectral convergence and sub-element resolution of the fields.
The DG approximation allows meshing of corners that are not multiples of $\pi/2$ in a discretization consistent manner, when needed.
The high order field can then be exploited to accurately find irregular {\color{black}nodes}, and can be accurately integrated using a high order {\color{black}separatrix} integration method to avoid features like limit cycles.
The result is a mesh with naturally curved quadrilateral elements that do not need to be curved \textit{a posteriori} to eliminate invalid elements.
The mesh generation procedure is implemented in the open source mesh generation program \emph{NekMesh}.

\end{abstract}


\end{frontmatter}

\section{Introduction}

A quadrilateral subdivision of a two dimensional domain into a minimal number of subdomains can serve as the starting point from which to generate meshes for finite element~\cite{Viertelabs-1708-02316}, finite volume, (block structured) finite difference~\cite{kozdon:2013} and spectral element methods~\cite{Karniadakis:2005fj,Kopriva:2009nx}.
To generate finite element and finite volume meshes, one would subdivide each subdomain into quadrilateral cells sufficiently small to resolve the geometry and expected solution features.
The resulting mesh would have the desired property of only a small number of irregular nodes where the valence, i.e., the number of elements sharing a point, is not equal to four.
The subdivision of the subdomains could also be used as nodes for high order block-structured finite difference methods if the subdivision within each block is structured.
Finally, the subdomains could be used \textit{as is} with high order geometry information intact, or subdivided in a structured or unstructured manner, to be used as elements for a quadrilateral based spectral element method.

\begin{figure}[htbp]
   \centering
   \includegraphics[width=0.49\textwidth]{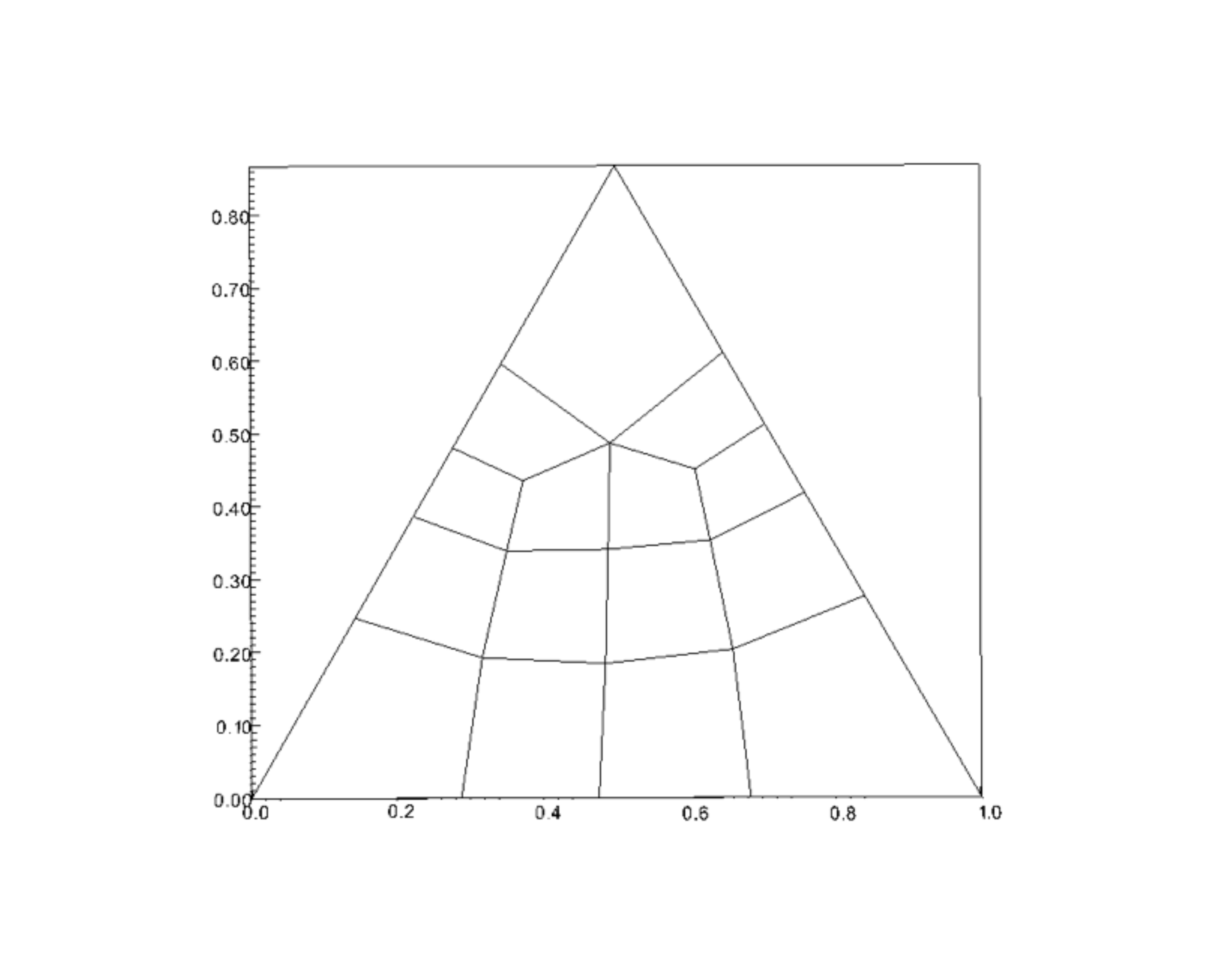} 
   \includegraphics[width=0.49\textwidth]{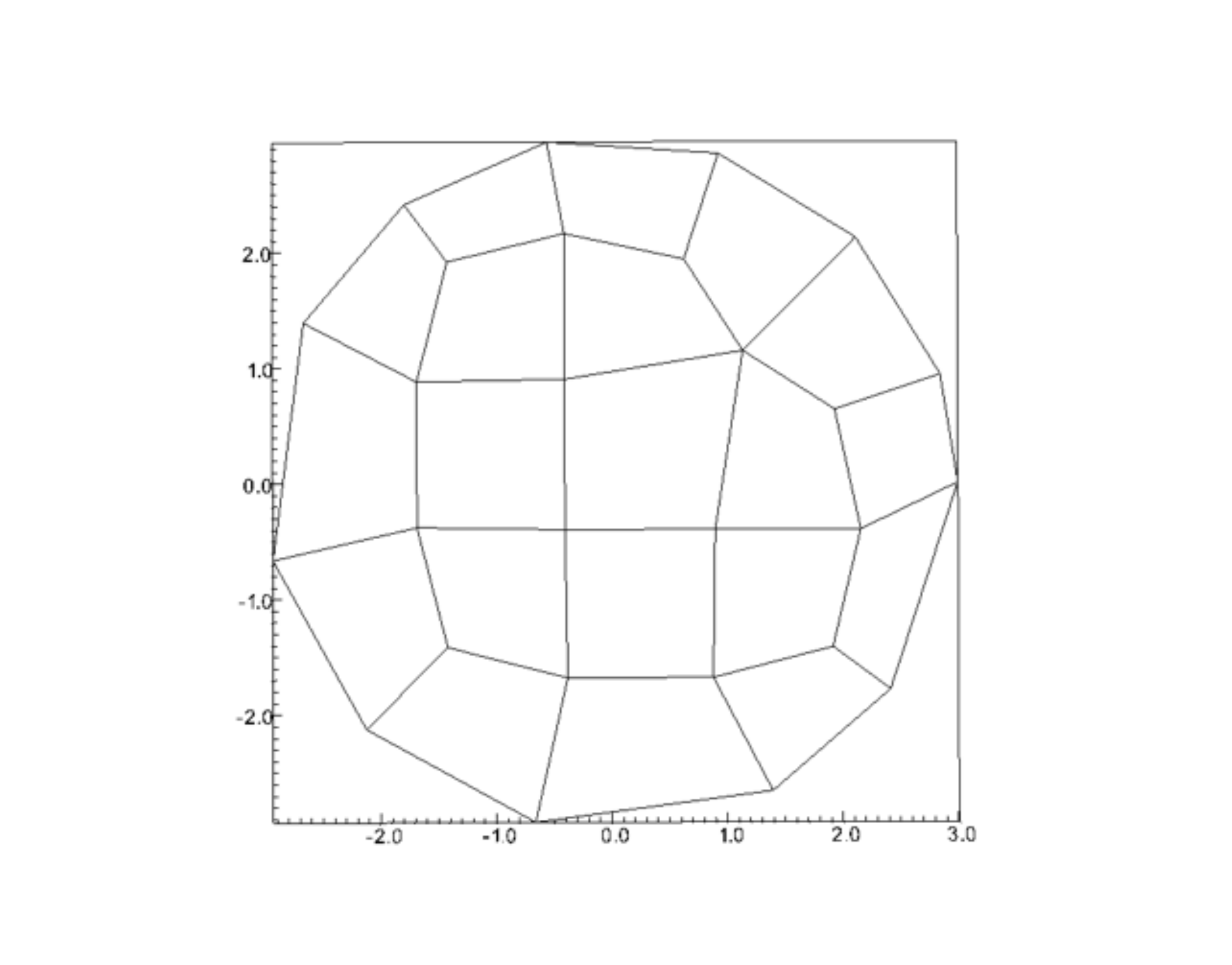} 
   \caption{Examples of quadrilateral meshes in simple geometries generated by a subdivision algorithm.}\label{fig:SubdivisionMeshes}
\end{figure}

Fully automatic quadrilateral mesh generation is still difficult~\cite{Bommes:2013pd}, especially for high order (curved) elements.
Most previous work has focussed on generating large numbers of small straight sided elements.
The extension to high order is usually done by curving the elements \textit{a posteriori}, e.g.~\cite{Hindenlang:2014gl}.
Unstructured quadrilateral mesh generators tend to generate meshes with large numbers of irregular {\color{black}nodes}~\cite{B.-H.-V.-Topping:2004aq}.
Extraneous irregular {\color{black}nodes} are generated even for simple geometries, where mesh topologies can be generated by hand with a minimum number of them.
As an example, Fig.~\ref{fig:SubdivisionMeshes} shows two meshes generated with a subdivision algorithm~\cite{liang2010,Schneiders-vki}, which generates more subdivisions than necessary.

\begin{figure}[htbp]
   \centering
   \includegraphics[width=0.49\textwidth]{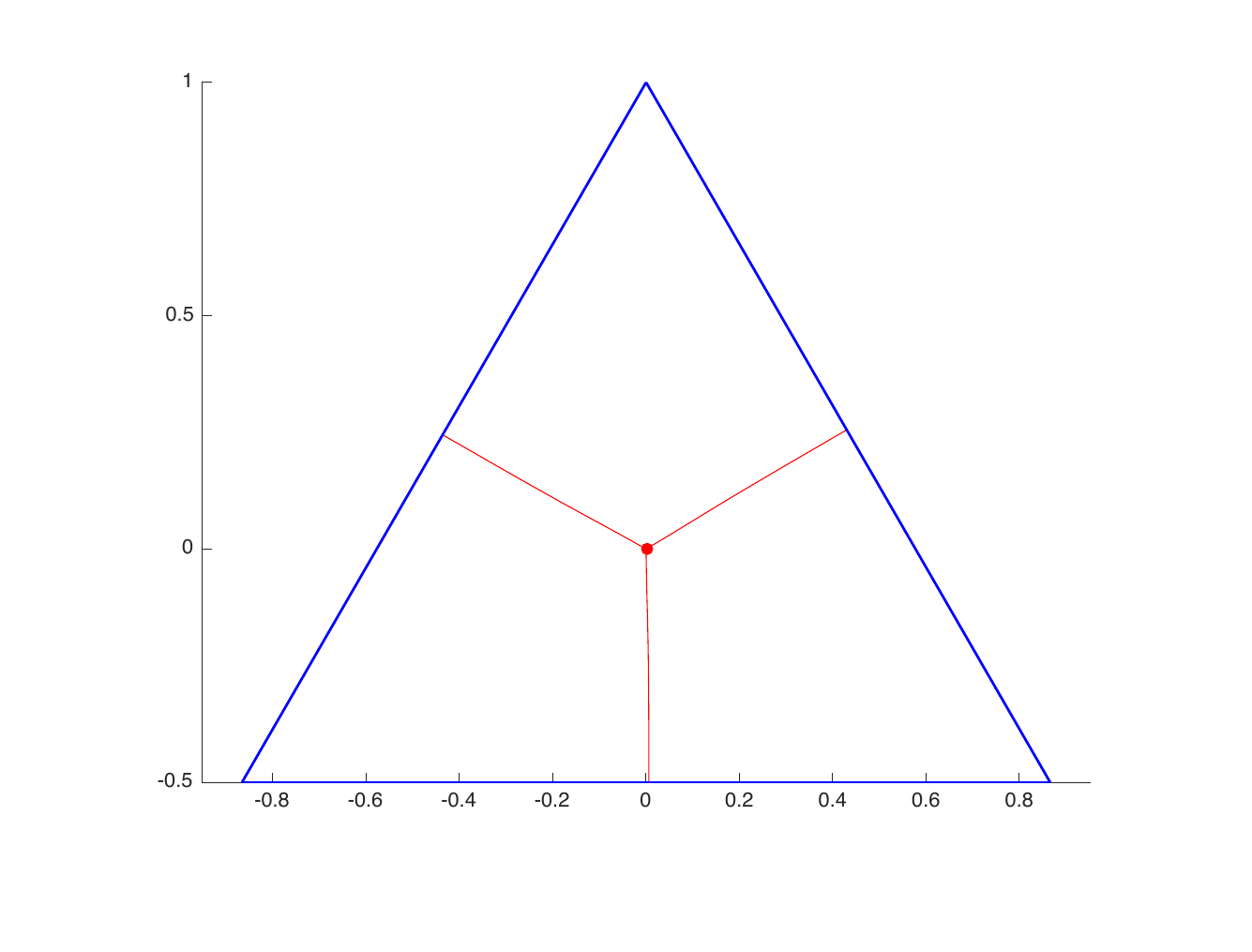} 
   \includegraphics[width=0.49\textwidth]{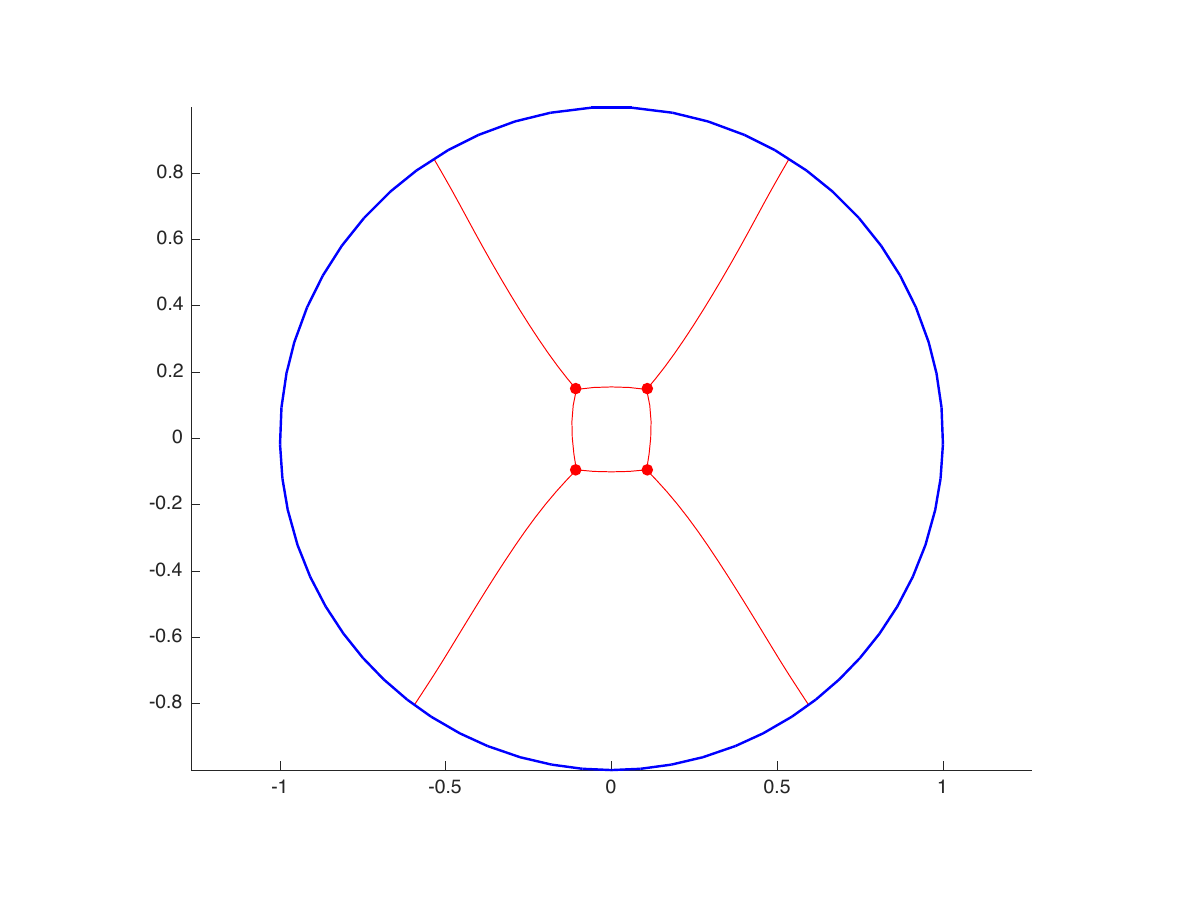} 
   \caption{Subdivision of the triangle and circle geometries using a cross field approach.}\label{fig:CrossFieldMeshes}
\end{figure}

One promising approach for generating quadrilateral subdivisions in arbitrary domains has come out of the computer graphics community.
It uses cross fields~\cite{Viertelabs-1708-02316,Bommes:2013pd} in a field-guided procedure.
Cross fields resemble cross hatching used in drawings.
Crosses, which are composed of two direction vectors and their negatives at a point, are invariant to rotations of $\pi/2$.
Cross field methods generate meshes with fewer irregular {\color{black}nodes}.
Fig.~\ref{fig:CrossFieldMeshes} shows the same two geometries as Fig.~\ref{fig:SubdivisionMeshes} decomposed with a cross field method.
The method produces the same topology as an experienced user might generate by hand, and has a minimum number of irregular {\color{black}nodes}.

Details of cross field procedures vary, but the basic idea is that the crosses are oriented, usually at the nodes of an existing triangular mesh, according to a minimization or smoothing procedure.
Triangles in which the crosses cannot smoothly vary contain singularities, which serve as the irregular nodes of the quadrilateral decomposition to be computed.
The internal edges of the quadrilaterals are found by integrating {\color{black}\emph{separatrices} (represented as \emph{streamlines} of the cross field)}~\cite{Nicolas-Kowalski:2012fu}, starting from the {\color{black}irregular nodes} until they reach another {\color{black}irregular nodes} or a physical boundary.
The quadrilateral subdivision is usually further subdivided into smaller quadrilateral elements.

Low order techniques {\color{black} using crosses} do not precisely locate {\color{black}irregular nodes} within the domain. {\color{black}Irregular nodes can be located only to within the size of the element in which they occur. 
Spurious singularities can be generated, which may have been coalesced in an \textit{ad hoc} fashion. For instance, if the underlying mesh is too coarse it is possible that two or more singularities can fall within one element.}
The valence of the {\color{black}irregular points} computed from a coalesced singularity could lead to the incorrect number of {\color{black}separatrices} and hence the incorrect valence of the associated interior quadrilateral mesh node.

Low order tracing of the {\color{black}separatrices} also often leads to multiple curves that need to be coalesced.
Associated with this is the well-known limit cycle problem where {\color{black}separatrices}, instead of meeting as they should, spiral indefinitely~\cite{Viertel:2017rt}.

Finally, boundary conditions are difficult to apply at corners where a jump discontinuity in the cross field will occur when a standard continuous Galerkin (finite element) approximation is used for smoothing.
Rather than putting a singularity there, the corner may be effectively smoothed by averaging the cross field at such points~\cite{Viertelabs-1708-02316}, at the expense of adding more interior singular points.

We address issues found in the use of low order techniques by using a high resolution approach to compute a guiding field, accurately locating the interior {\color{black}irregular nodes}, and accurately integrating the streamlines used to subdivide the domain.
We use a high order continuous or discontinuous Galerkin spectral element method (SEM or DGSEM) on a triangular mesh, the approximation depending on the regularity of the boundary conditions needed to solve Laplace problems.
Unlike previous work, we perform all operations on the original highly-resolved guiding field rather than a cross field so that accuracy is not lost.
We use an accurate locator for the {\color{black}irregular} points that exploits the high resolution solution, and an accurate method for computing their valence.
The streamline integration then uses the high order information available to obtain precise streamlines, which reduces the need to coalesce extraneous lines.

\section{Mathematical Formulation}\label{sec:MathematicalFormulation}

\begin{figure}[htbp]
   \centering
   \includegraphics[width=2.15in]{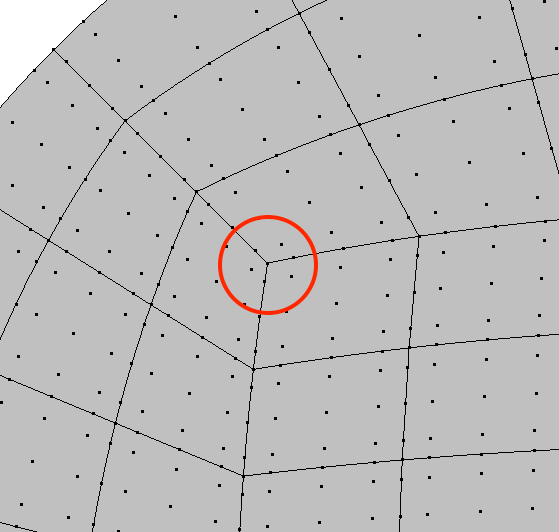} 
    \hspace{0.5cm}
   \includegraphics[width=2in]{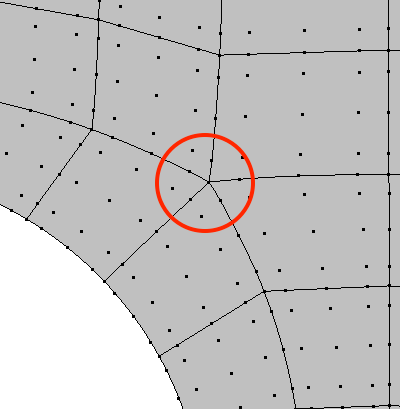} 
   \caption{Portions of quadrilateral meshes, where irregular nodes are encircled with red. Left: Valence three. Right: Valence five. Interior degrees of freedom for fourth order spectral elements are represented by dots.}\label{fig:ValenceExample}
\end{figure}

The problem is to subdivide any two dimensional domain, $\Omega$, that has piecewise smooth boundaries into quadrilateral subdomains.
The domain might be simply or multiply connected.
The decomposition will be regular if all the corners of the quadrilaterals have valence four, that is, each corner node connects four edges.
In general, it will not be possible to generate a purely regular mesh, especially in multiply connected domains.
Instead, some of the nodes will be \emph{irregular}, where the valence will be larger or smaller than four, like those observed in Fig.~\ref{fig:ValenceExample}.

\subsection{Crosses and Cross Field}

\begin{figure}[htbp]
   \centering
   \includegraphics[width=3in]{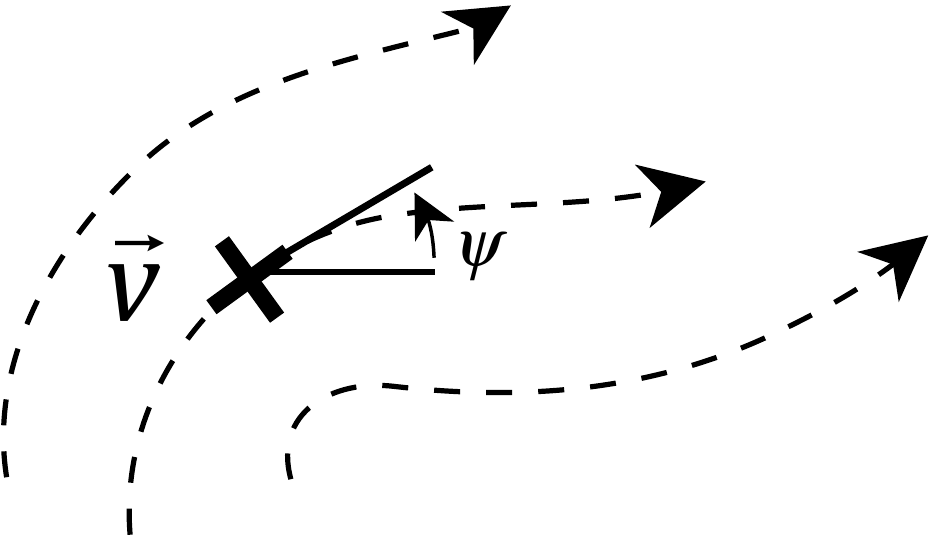} 
   \caption{A cross defined in a field $\vec v$. The axes of the cross lie in the tangent and normal directions to the field lines.}\label{fig:CrossOnField}
\end{figure}

To motivate the high resolution approach to the quadrilateral decomposition of a domain, $\Omega$, and to put it into context, we recall the notion of a cross.
A cross can be represented by
{\color{black}
\begin{equation}
C\left( \psi  \right) = \left\{ {{{\vec c}_k} = {{\left( { c_{ x}, c_{y}} \right)}_k} = \left( {\cos \left( {\psi  + k \frac{\pi}{2}} \right),\sin \left( {\psi  + k \frac{\pi}{2}} \right)} \right)} \right\},
\label{eq:CrossDef}
\end{equation}
for, $k = 0,1,2,3$.}
Here, the tangent angle, or phase, $\psi$, is computed from the four quadrant inverse tangent,
\begin{equation}
\psi = \frac{1}{4}\rm{atan2}\left({v,u}\right)\in\left[-\frac{\pi}{4},\frac{\pi}{4}\right],
\label{eq:PsiDef}
\end{equation}
so that a cross is represented by four unit vectors {\color{black}$\vec c_{k}$} at any point $\vec x \in \Omega$ given a vector field $\vec v = \left(u\left(\vec x\right),v\left(\vec x\right)\right)$ at which $\vec v\ne 0$ from the principal direction.
See Fig.~\ref{fig:CrossOnField}.
The important property of a cross is its 4-way rotational symmetry.
Rotation by an angle $\pi/2$ does not change the cross.

The function $\psi$ is tangent or orthogonal (due to the jump in the atan2 function) to the streamlines of the field $\vec v$.
Hence, it is parallel to one of the branches of a cross.
It is undefined at $\vec v = (0,0)$.
Since $\psi$ is computed from the arctangent, there will be a jump of of value $\pi/2$ in $\psi$ depending on the signs of $u$ and $v$ even if $\vec v$ is smooth.
We call lines in the field across which $\psi$ jumps \emph{jump lines}.

A full discussion of cross fields and methods based on them is beyond the scope of this paper.
See, for example, references~\cite{Viertelabs-1708-02316,Viertel:2017rt,Bunin2006,2008arXiv0802.2399B} for more comprehensive discussions.

\subsection{The Guiding Field, $\vec v$}

\begin{figure}[htbp]
  \centering
  \includegraphics[width=2in]{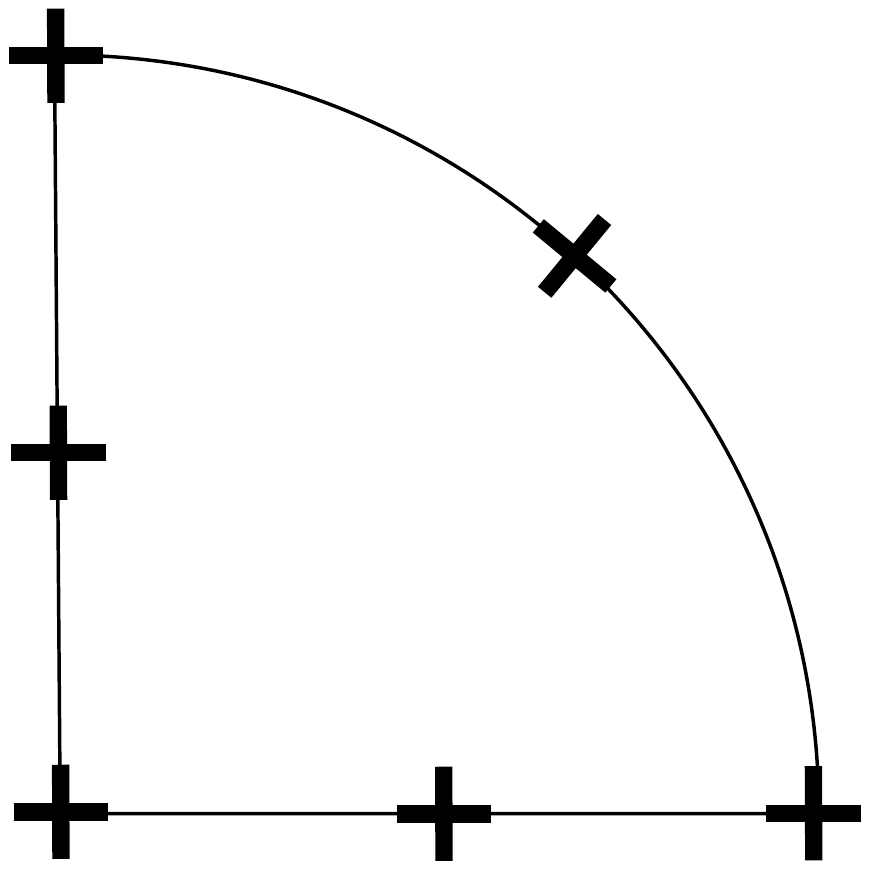} 
  \caption{Boundary crosses shown on a quarter circle domain. Along the circular boundary portion, $u_{b}$ has two zeros and one zero contour, while  $v_{b}$ has three zeros.}\label{fig:CrossesOnQuarter.pdf}
\end{figure}

It is necessary, then, to find the guiding field $\vec v$ from which to find the {\color{black}irregular nodes} and trace {\color{black}separatrices}.
The only values that can be set \textit{a priori} are on the boundary $\partial \Omega$ of $\Omega$.
We align the field with the boundary to ensure that the mesh is aligned there.
Therefore, we set
\begin{equation}
{{\vec v}_b} = \left(u_{b},v_{b}\right)=\left( {\cos \left(  4\theta_{b} \right),\sin ( 4\theta_{b})} \right),
\label{eq:BCs}
\end{equation}
where $\theta_{b}$ is the tangent angle of the boundary where the field is being computed.
The vector $\vec v_{b}$ at any point along the boundary defines a cross~\eqref{eq:CrossDef} at that point.
The factor of four in eq.~\eqref{eq:BCs} ensures the same $u,v$ values for each $90^{\circ}$ rotation of the angle, and hence the rotational symmetry.
Fig.~\ref{fig:CrossesOnQuarter.pdf} shows example crosses on the boundary of a quarter circle domain.

We make three observations from eq.~\eqref{eq:BCs} about the boundary field $\vec v_{b}$:  
\begin{itemize}
  \item Along any smooth closed boundary curve (interior or exterior), there are eight zeros each for $u$ and $v$, and those zeros are not coincident.
  \item Along smooth portions of the boundary curves and at corners whose angle is a multiple of $\pi/2$, $u$ and $v$ are continuous.
  \item At boundary corners with angles not divisible by $\pi/2$ there is a jump discontinuity in $\vec v_{b}$ and the cross field.
\end{itemize}

Once the boundary values are specified using expression~\eqref{eq:BCs}, they are smoothly propagated to the interior.
Smoothing procedures that have been used in the past include marching the associated crosses to the interior~\cite{Fogg2015}, minimizing an energy functional~\cite{Viertelabs-1708-02316}, or solving a Laplace problem~\cite{J.-F.-Remacle1:2013ca}.
Solving a Laplace problem for $u$ and $v$ 
\begin{equation}
\left\{ \begin{gathered}
  {\nabla ^2}\vec v = 0,\quad \, \vec x  \in \Omega  \hfill \\
  \vec v = {{\vec v}_b},\quad \, \vec x  \in \partial \Omega  \hfill \\ 
\end{gathered}  \right.
\label{eq:LaplaceEqns}
\end{equation}
guarantees that the field $\vec v$ is smooth in the interior of the domain and satisfies the maximum principle.
The solution of the Laplace problem is equivalent to the minimization of the Dirichlet energy.
We choose to solve~\eqref{eq:LaplaceEqns} because it allows us to use a high resolution approach such as the spectral element methods.
In particular, the use of a discontinuous Galerkin approximation allows us to handle geometries where the corner angles generate {\color{black}jumps} in the guiding field, see Sec.~\ref{sec:CGAndDG}.

\subsubsection{Continuous and Discontinuous Galerkin Spectral Element Methods}\label{sec:CGAndDG}

A key difference between our work here and that of previous work is that we compute high resolution solutions for the BVPs~\eqref{eq:LaplaceEqns} with either a continuous Galerkin (CG) or a discontinuous Galerkin (DG) spectral element method on a triangular mesh~\cite{Karniadakis:2005fj}.
Spectral element methods are spectrally accurate, meaning that the convergence rate depends only on the smoothness of the solution.
They are high resolution in that they use a large number of degrees of freedom within an element.
Unlike traditional low order finite element methods, the high order polynomial expansion of the solution inside each element allows us to locate {\color{black}irregular nodes}, identify their valence, and finally trace {\color{black}separatrices} with high order accuracy.

In this work, we use the spectral/\textit{hp} element methods formulation described in detail in reference~\cite{Karniadakis:2005fj} and implemented in \emph{Nektar++}~\cite{Cantwell2015,Moxey2019}.
We briefly describe the fundamentals of the methods in what follows.
We consider the numerical solution of PDEs of the form \(\mathcal{L} u = 0\) over a domain \( \Omega \).
The domain \(\Omega \) is taken as a set of finite elements, \(\Omega_e\) --- the mesh --- such that \(\Omega = \cup \Omega_e\) and \(\Omega_{e_1} \cap \Omega_{e_2} = \partial\Omega_{e_1 e_2}\) is either an empty set or the interface between two elements and is of one dimension less than the mesh.
We solve the PDE problem in the weak sense and require that \(\left.u\right|_{\Omega_e}\) is in the Sobolev space \(H^{1}(\Omega_e)\).
In a CG formulation, we also require the solution to be continuous from one element to the other.

We formulate the solution to~\eqref{eq:LaplaceEqns} in weak form: \(find\ \vec v \in H^1 (\Omega)\ such\ that\)
\begin{equation}
  a(\vec v, \vec w) = l(\vec w) \quad \forall \vec w \in H^1 (\Omega)
\end{equation}
where \(a(\cdot,\cdot)\) is a symmetric bilinear form, \(l(\cdot)\) is a linear form and \(H^1 (\Omega)\) is formally defined as
\begin{equation}
  H^1 (\Omega) = \{ \vec w \in L^2 (\Omega) \mid D^\alpha \vec w \in L^2 (\Omega)\ \forall\ \left|\alpha\right| \leq 1 \}.
\end{equation}

We want to solve this problem numerically and therefore consider solutions in a finite dimensional subspace \(V_N \subset H^1(\Omega)\) and cast the problem: \(find\ \vec v^\delta \in V_N\ such\ that\)
\begin{equation}
  a(\vec v^\delta, \vec w^\delta) = l(\vec w^\delta) \quad \forall \vec w^\delta \in V_N,
\end{equation}
augmented with boundary conditions,~\eqref{eq:BCs}.
In the CG formulation, we also enforce the condition \(V_N \subset C^0\).

We take a weighted sum of \(N\) trial functions \(\Phi_n( \vec x )\) defined on \(\Omega \) so we can represent \(\vec v^\delta( \vec x ) = \sum_n \hat{v}_n \Phi_n( \vec x )\).
This transforms the problem to that of finding the coefficients \(\hat{u}_n\) that define \(\vec v^\delta( \vec x )\) within an element.
To obtain a unique choice of coefficients \(\hat{u}_n\), we place a restriction on \(R = \mathcal{L} \vec v^\delta \) that its \(L^2\) inner product, with respect to the test functions \(\Psi_n( \vec x )\), is zero.
In the Galerkin projection one chooses the test functions to be the same as the trial functions, i.e. \(\Psi_n = \Phi_n\).

The contributions of each element in the domain must be taken into account to construct the global basis \(\Phi_n\).
A parametric mapping \(\mathcal{X}_e : \Omega_e \to \mathcal{E}\) exists from each element \(\Omega_e\) to a standard reference element \(\mathcal{E} \subset {[-1,1]}^d\).
This mapping is given by \( \vec x  = \vec{\mathcal{X}}_e\left( \vec \xi\right)\).
It is important to distinguish \( \vec x \), the physical coordinates, from the \( \vec \xi \), the coordinates in the reference space.

We construct a local polynomial basis on the reference element to represent the solution.
A one-dimensional order-\(P\) basis is a set of polynomials \(\Phi_p(\xi), 0 \leq p \leq P\), defined on the reference segment \(-1 \le \xi_1 \le 1\).
In two and three dimensional reference regions, a tensor basis is used, where the polynomial space is constructed as the tensor product of one dimensional bases on segments, quadrilaterals, and hexahedral reference regions.
{\color{black}
In \emph{Nektar++}, triangular, tetrahedral, prismatic, and pyramidal elements are created by collapsing one or more of the coordinate directions to create singular vertices.
This allows it to support, for this work, easy-to-generate triangular meshes.
}

Finally, the discrete solution in a physical element \(\Omega_e\) can be expressed as
\begin{equation}
\vec v  ^\delta( \vec x ) = \sum_n \hat{v}_n \phi_n (\mathcal{X}_e^{-1} ( \vec x )),
\label{eq:SolutionInterpolant}
\end{equation}
with \(\hat{v}_n\) the coefficients computed by the Galerkin procedure.
We restrict our solution space to
\begin{equation}
  V := \{\vec v \in H^1(\Omega) \mid \left.\vec v\right|_{\Omega_e} \in \mathcal{P}_P(\Omega_e)\},
\end{equation}
where \(\mathcal{P}_P(\Omega_e)\) is the space of order \(P\) polynomials on \(\Omega_e\).

An assembly operator can then be designed to assemble the element contributions to the global solution.
In the CG formulation, elemental contributions of neighbours are summed to enforce \(C^0\)-continuity.
In the DG formulation, such mappings transfer flux values from the element interfaces into the global solution vector.

\subsection{{\color{black}Critical Points}}\label{sec:math_critical}

Points where $v = u = 0$ are called \emph{critical points} in the guiding field and \emph{singular points} in the cross field.
They are points where $\psi$ cannot be determined uniquely.
Such points become the interior \emph{irregular nodes} (valence $\ne 4$) in the quadrilateral decomposition.

\begin{figure}[htbp]
  \centering
  \includegraphics[width=0.49\textwidth]{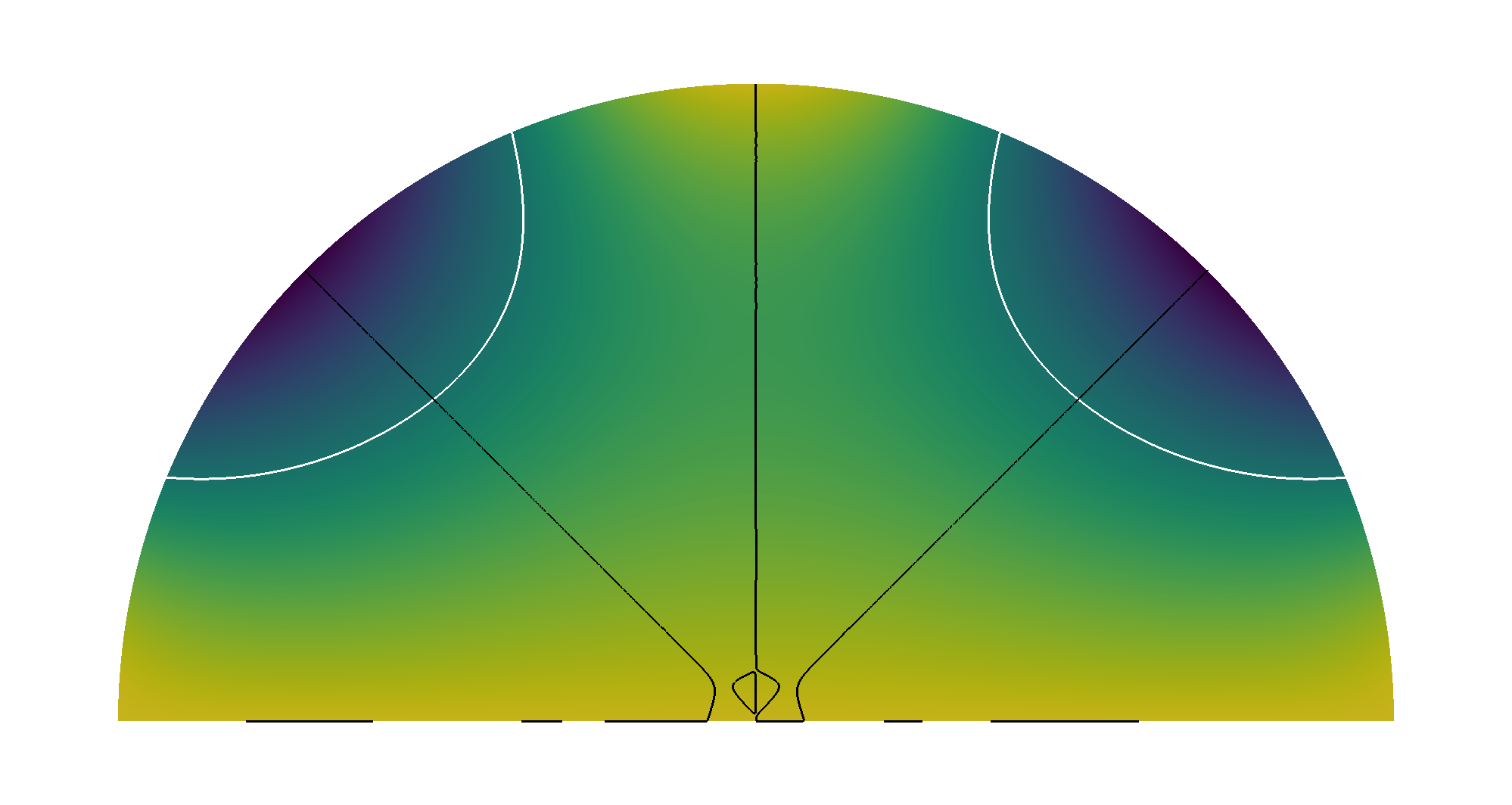}
  \includegraphics[width=0.49\textwidth]{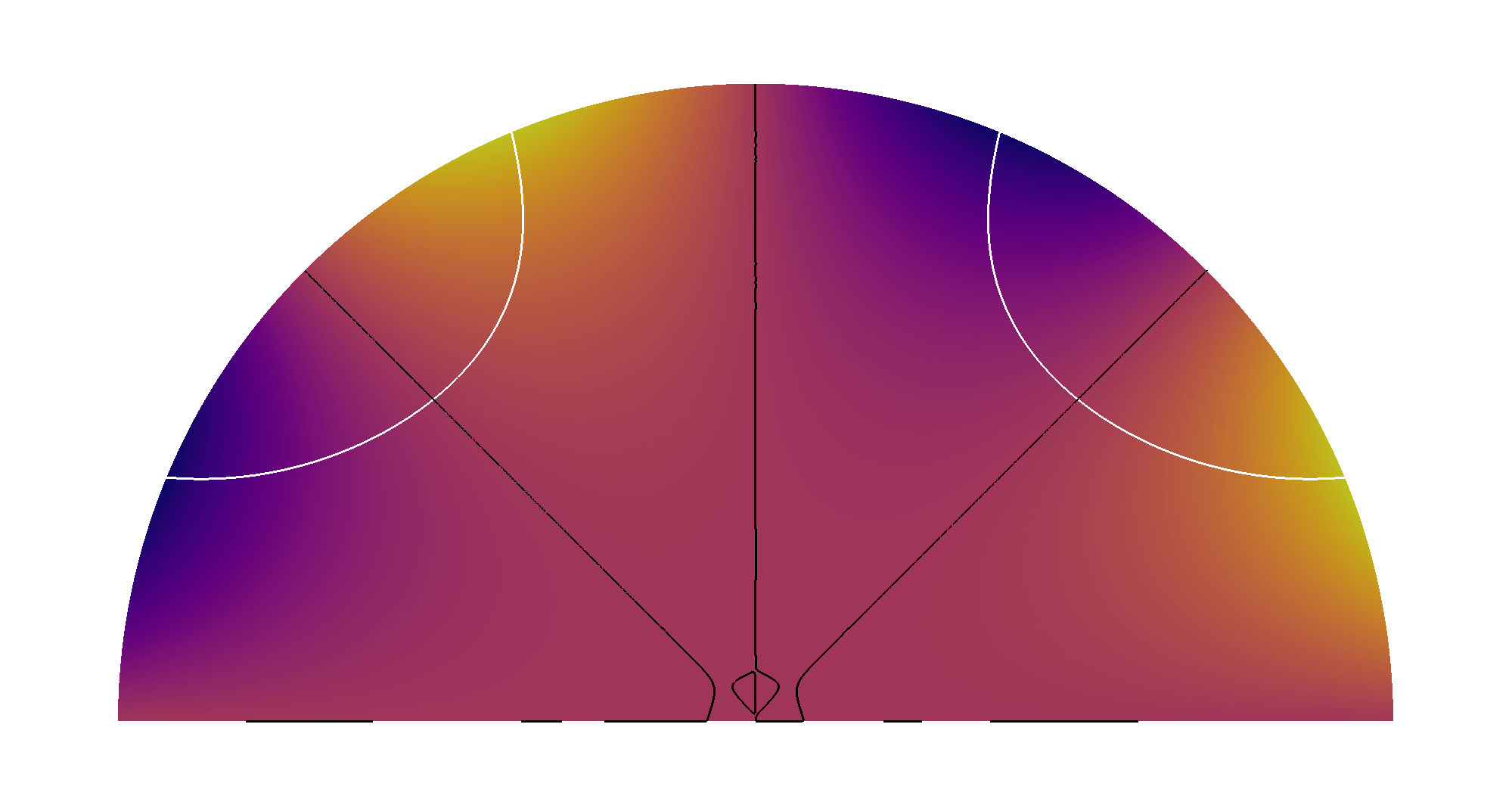}
  \includegraphics[width=0.75\textwidth]{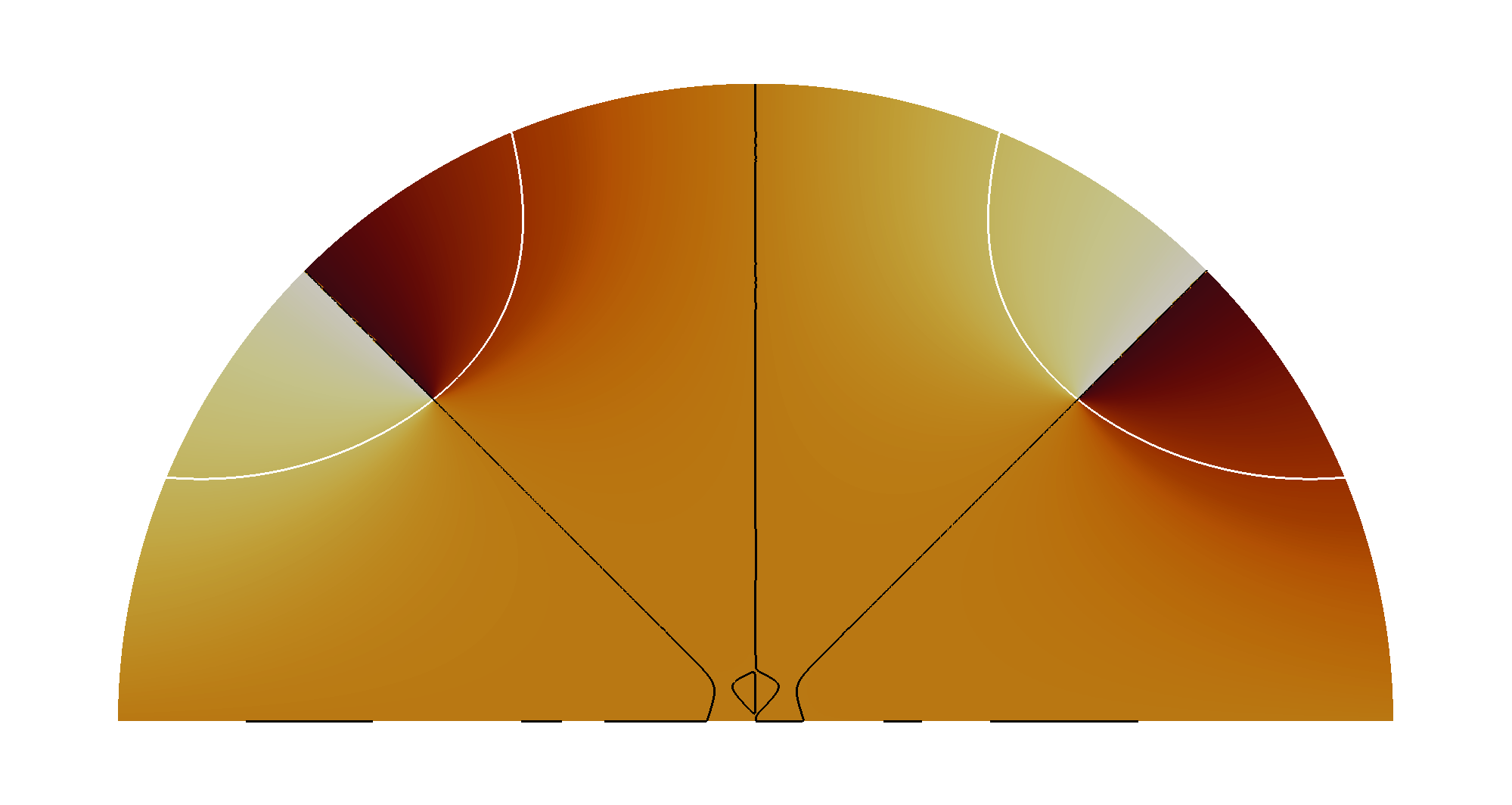}
  \caption{Solution \(u\) (top left) and \(v\) (top right) fields and computed \(\psi \) field (bottom) on the half disc geometry. Isocontours of \(u = 0\) and \(v = 0\) shown in white and black respectively.}\label{fig:half-disc-solution}
\end{figure}

The form of the boundary values usually implies that {\color{black}critical} points will exist in the interior of the domain.
{\color{black}Since} there are {\color{black}interlaced} zeros in $u$ and $v$ on the boundary, such as occur when the boundary is smooth and $\theta_{b}$ varies continuously, the fact that the interior field is smooth implies that zero contours of $u$ and $v$ must cross somewhere in the interior, thereby creating a {\color{black}critical} point.
Zeros in $u$ will exist along the boundary at any point where $\theta_{b}$ smoothly passes through $\theta_{b}=k\pi/8$, $k$ odd, and in $v$ for $k$ even.
Therefore, one can infer the existence of interior irregular {\color{black}nodes} in the final quadrilateral mesh from the boundary curves.
A semi-circular section of a boundary, for instance will have four zeros in $u_{b}$ and five in $v_{b}$ and, in the absence of other nearby features, will create two interior {\color{black}critical} points.
See Fig.~\ref{fig:half-disc-solution}.
We note that although the cross field directions are undefined at {\color{black}critical points}, the field $\vec v$ is smooth because the solutions of the Laplace equation~\eqref{eq:LaplaceEqns} are regular.

{\color{black}Irregular nodes in the mesh} occur as recognizable {\color{black}critical point} topologies in the $\vec v$ field~\cite{Tricoche:2001:CTS:601671.601695}.
They can be categorized analytically by the Poincar\'e index
\begin{equation}
  {i_\gamma } = \frac{1}{{2\pi }}\oint_\gamma  {d\phi },
\end{equation}
where
{\color{black}
\begin{equation}
  \phi = \rm{atan2}\left({v,u}\right) = 4 \psi.
  \label{eq:phi-psi}
\end{equation}
}

When $\gamma$ encloses a {\color{black}critical} point, $i_{\gamma}=\pm 1$, where $-1$ corresponds to a saddle point.
Other critical points have index $+1$.
If there is no critical point inside the contour $\gamma$ then the index is zero.
The categorization of the critical point lets one determine the number of streamlines to emanate from the {\color{black}irregular node}, and hence its valence in the mesh.
When $i_{\gamma} = -1$ there are five separatices.
When $i_{\gamma} = +1$ there are three~\cite[Lemma 5.1]{Viertelabs-1708-02316},~\cite{Nicolas-Kowalski:2012fu}. {\color{black}It is possible to have an index of $+2$ at which six elements would meet. However such points are unstable and in practice will split into separate critical points
of lower valence. For this reason, we do not consider such critical points here.}

Using the definition~\eqref{eq:PsiDef} and relation~\eqref{eq:phi-psi}, we similarly define the integral
\begin{equation}
  I_{c} = \frac{1}{\pi/2} \oint_c  {\frac{d\psi}{d\theta}d\theta },
  \label{eq:IDef}
\end{equation}
where $c$ is a (small) counter-clockwise circular contour centered on the {\color{black}critical} point.
Table~\ref{tab:IValues} shows the values of $I_{c}$ for three {\color{black}separatrix} topologies in the neighborhood of a {\color{black}critical point}.
It is equal to the sum of the jumps in $\psi$ around the contour.
The \emph{valence} \(\cal V\) of the critical point is then
{\color{black}
\begin{equation}
  \mathcal V= 4 - I_{c}.
  \label{eq:critical_valence}
\end{equation}
}
Note, for consistency, that if there is no critical point inside the contour then $I_{c}=0$ and the point is regular with valence four.

\begin{table}[htbp]
  {\color{black}
  \begin{center}
    \caption{Value of $I_{c}$,~\eqref{eq:IDef}, and the associated valence.}\label{tab:IValues}

    \begin{tabular}{cc}$I_{c}$ &Valence ($\mathcal V$) \\
    \hline
    $-1$ & 5 \\
    $0$ & 4 \\
    $+1$ & 3
    \end{tabular} 
  \end{center}
  }
\end{table}
 
If one generalizes the integral $I_{c}$ to be over a part of the circle, one can characterize boundary vertices, too~\cite{2008arXiv0802.2399B}.
Let
\begin{equation}
  {I\left(\theta_{0},\theta_{f}\right) } = \frac{1}{\pi/2} \int_{{\theta _0}}^{{\theta _f}} {\frac{{d\psi }}{{d\theta }}d\theta },
\end{equation}
so that $I_{c}=I\left(0,2\pi\right)$.
At a boundary point, then, the number of quadrilaterals attached to the point (the boundary itself already being a streamline by construction) is
\begin{equation}
  \mathcal V = {\frac{\Delta\theta}{\pi/2} - {I\left(\theta_{0},\theta_{f}\right) } } ,
  \label{eq:math_corner}
\end{equation}
where \( \Delta\theta = \theta_f - \theta_0 \).
\( I\left(\theta_{0},\theta_{f}\right) \) can be viewed as a correcting factor to \( \Delta\theta \) towards a multiple of \( \frac{\pi}{2} \).

\subsection{Streamline Integration}

The guiding field is also used to trace \emph{streamlines} that will form the \emph{separatrices} of the block decomposition.
We can formulate the problem as finding the trajectory \( \vec x (t)\) that satisfies
\begin{equation}
  \frac{d \vec x }{dt} = \tilde v\left(\psi'( \vec x )\right)
  \label{eq:streamline-integration}
\end{equation}
where \(t\) is the integration parameter and \(\tilde v\left(\psi'( \vec x )\right)\) is the adjusted guiding field vector.
Equation~\eqref{eq:PsiDef} gives a guiding field \(\psi \in \left[-\frac{\pi}{4},\frac{\pi}{4}\right]\), but crosses are invariant to rotations of \(\pi/2\).
It is therefore necessary to account for the branches of the cross that don't lie within \(\left[-\frac{\pi}{4},\frac{\pi}{4}\right]\).
These adjusted guiding directions can be obtained by adding a certain number \(k = 0,1,2,3\) of \(\pi/2\) so that \(\psi'\) lies within the appropriate quadrant \(\left[-\frac{\pi}{4}+k\frac{\pi}{2},\frac{\pi}{4}+k\frac{\pi}{2}\right]\).
This adjustment operation is especially important when streamlines cross {jump lines} where \(\psi \) abruptly changes by \(\pm\pi/2\), but the streamline has to continue in the same overall direction.

The ordinary differential equation~\eqref{eq:streamline-integration} requires initial conditions in the form of
\begin{equation}
  \vec x (t_0) =  \vec x _0
\end{equation}
where \(t_0\) is the initial integration time and \( \vec x _0\) the location of the start of the streamline.
In this case, the start of the streamline will either be an {\color{black}irregular node} or a corner from which a non-boundary separatrix emerges.
The ODE~\eqref{eq:streamline-integration} can then be integrated in a traditional numerical manner as will be described below in Section~\ref{sec:streamlines}.

\section{The Mesh Generation Process}\label{sec:implementation}

From the discussion above, we see that the decomposition of a domain into quadrilateral subdomains has four stages:
\begin{enumerate}
  \item Computation of a guiding field.
    This includes the generation of a triangular mesh and solution of the Laplace problem~\eqref{eq:LaplaceEqns}.
  \item Finding {\color{black}critical points of the guiding field} and their valences, including the valences of corners at boundary vertices.
  \item Integration of streamlines, from {\color{black}critical points} and corners, to generate the interior boundaries of the quadrilateral decomposition.
  \item Cutting of the domain into quadrilaterals using the {\color{black}separatrices} as subdomain boundaries, then further subdivision into elements, as necessary.
\end{enumerate}

\begin{figure}[htbp]
  \centering
  \includegraphics[width=0.75\textwidth]{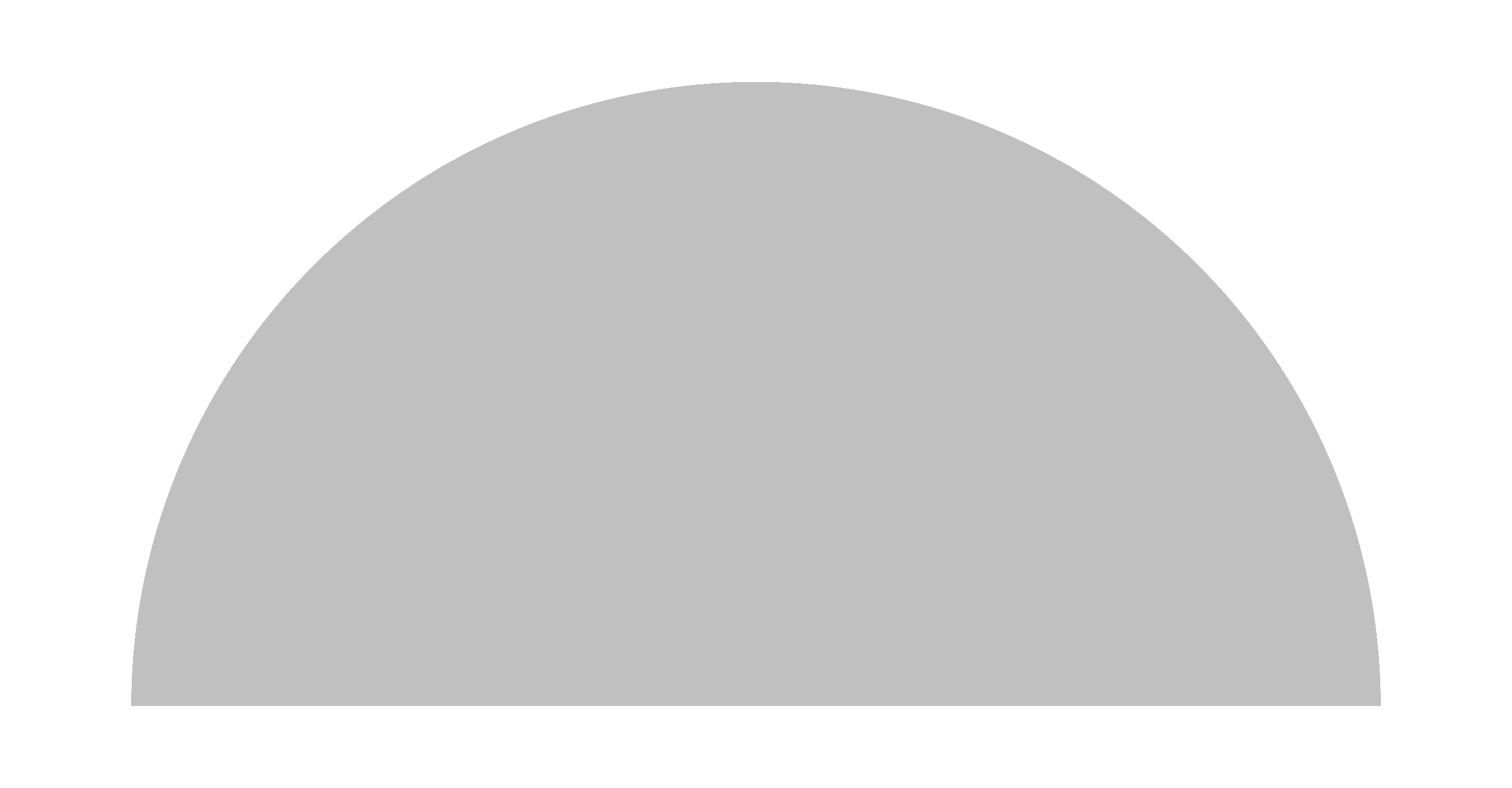}
  \caption{Geometry of a half disc for illustration of the quadrilateral meshing method.}\label{fig:half-disc-geometry}
\end{figure}

This section presents details of the implementation of the process and walks through the four stages needed to generate a quadrilateral mesh.
The procedure is illustrated using a geometry commonly used in the cross field literature~\cite{Nicolas-Kowalski:2012fu,Viertel:2017rt}, that of a half disc as shown in Fig.~\ref{fig:half-disc-geometry}.

\subsection{Solution of the Field Equations}\label{sec:solution}

\begin{figure}[htbp]
  \centering
  \includegraphics[width=0.75\textwidth]{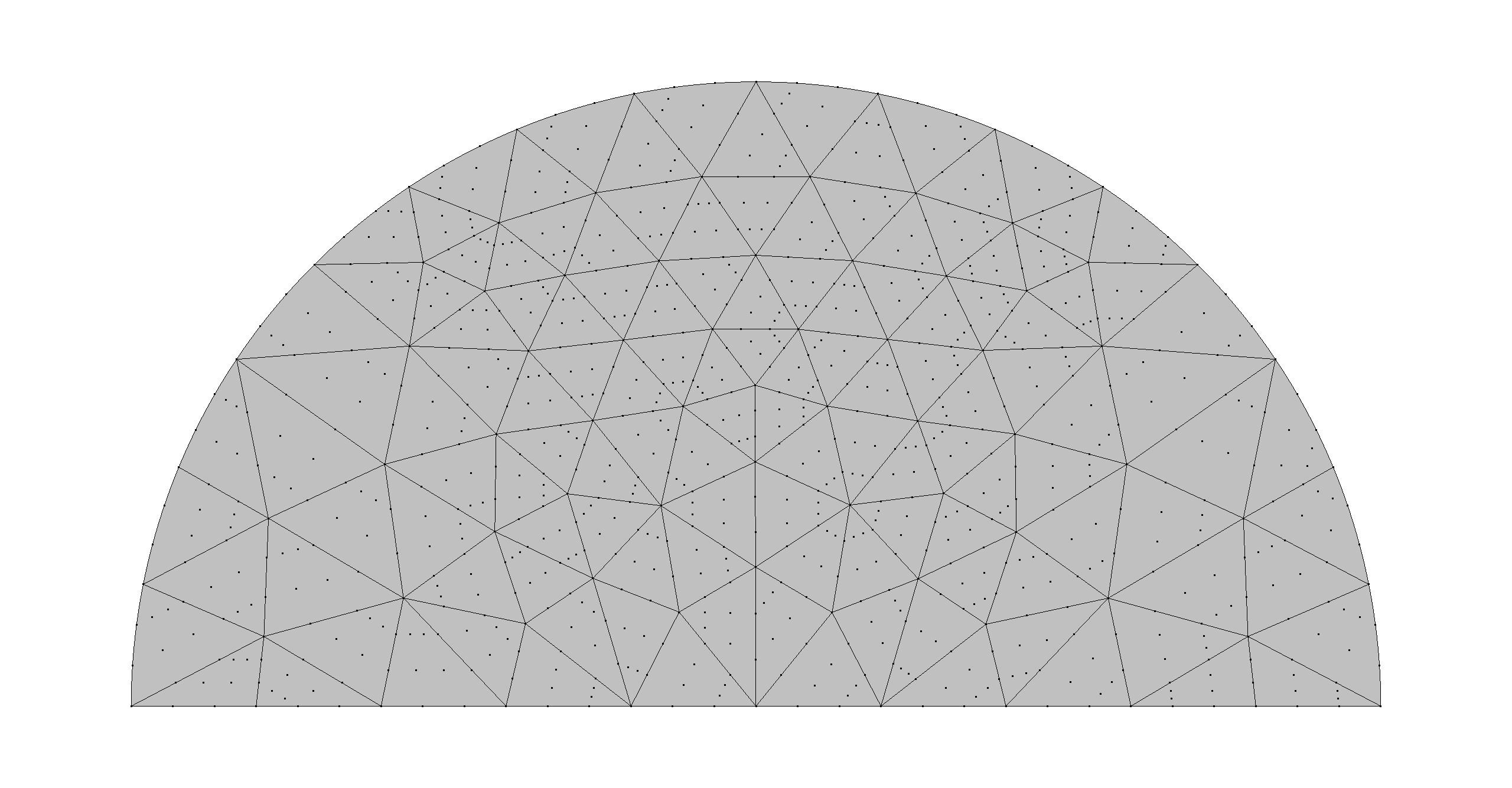}
  \caption{Third order triangular spectral element mesh for the half disc geometry.}\label{fig:half-disc-tris}
\end{figure}

To solve the Laplace problem~\eqref{eq:LaplaceEqns} we use the open source spectral element program \emph{Nektar++}~\cite{Cantwell2015,Moxey2019}.
We first generate a triangular finite element mesh which is then made high order and curved by projecting interior nodes onto the curved boundaries~\cite{Sherwin2002}.
This procedure is carried out in \emph{NekMesh}~\cite{Moxey2019}, which also has the capability to optimize the high order mesh should some elements be of low quality or simply invalid~\cite{Turner2017}.
Such a mesh is shown in Fig.~\ref{fig:half-disc-tris} for the reference geometry.

For the purposes of generating the quadrilateral decomposition, we generate as coarse a triangular mesh as possible and then use high order polynomials to approximate the solution within the elements to get the desired accuracy.
First, we want to take advantage of the spectral accuracy of the spectral element Laplace solver, where a fine mesh is not necessary.
Second, as will be discussed in Sections~\ref{sec:critical} and~\ref{sec:streamlines}, a small number of elements is desirable to simplify and speed up the {\color{black}critical point} detection and streamline tracing.
{\color{black}For the half circle geometry and triangle mesh shown in Fig.~\ref{fig:half-disc-tris}, it was sufficient to use third order polynomials (fourth order convergence with element size).}

Once the curvilinear triangular mesh is generated, the Laplace problem~\eqref{eq:LaplaceEqns} is solved.
A special implementation of the Laplace solver has been implemented in the spectral element code \emph{Nektar++}~\cite{Cantwell2015} for the boundary value problem~\eqref{eq:LaplaceEqns}, where the boundary conditions are automatically computed using~\eqref{eq:BCs}.

{\color{black}
\subsubsection{Choice of a Discretization}
}

\emph{Nektar++} supports both continuous Galerkin and discontinuous Galerkin discretizations.
The CG approximation is used when the boundary conditions (BCs) are continuous along all boundary curves.
As can be seen in the discussion in Sec.~\ref{sec:MathematicalFormulation},  boundary conditions are continuous if the boundaries satisfy one of two conditions:
\begin{itemize}
  \item The curve at each point is \(C^1\)-continuous, i.e.~it is a smooth curve;
    or 
  \item The boundary curve is only \(C^0\)-continuous, i.e.~a corner, and the angle is a multiple of \(\pi/2\).
\end{itemize}

In this paper, we employ CG whenever possible and switch to DG only when the geometry requires discontinuous BCs.
For example, the reference geometry of Fig.~\ref{fig:half-disc-geometry} contains only smooth curves and \(\pi/2\) corners.
We therefore use a CG formulation for it.
The solution for \(\vec v\) on the half disc geometry is shown in Fig.~\ref{fig:half-disc-solution}.
As expected, boundary conditions and the solution are smooth, so the solution benefits from rapid convergence of the SEM.\@

If neither smoothness condition is satisfied at a boundary corner point, the boundary condition~\eqref{eq:BCs} has a jump discontinuity there.
In that case, we use a DG formulation, which can account for discontinuous BCs.
Unlike the CG approximation used in traditional low order solvers, no \textit{ad hoc} smoothing of a corner's BCs is required, leading to a discretization-consistent solution.

\begin{figure}[htbp]
  \centering
  \includegraphics[width=0.49\textwidth]{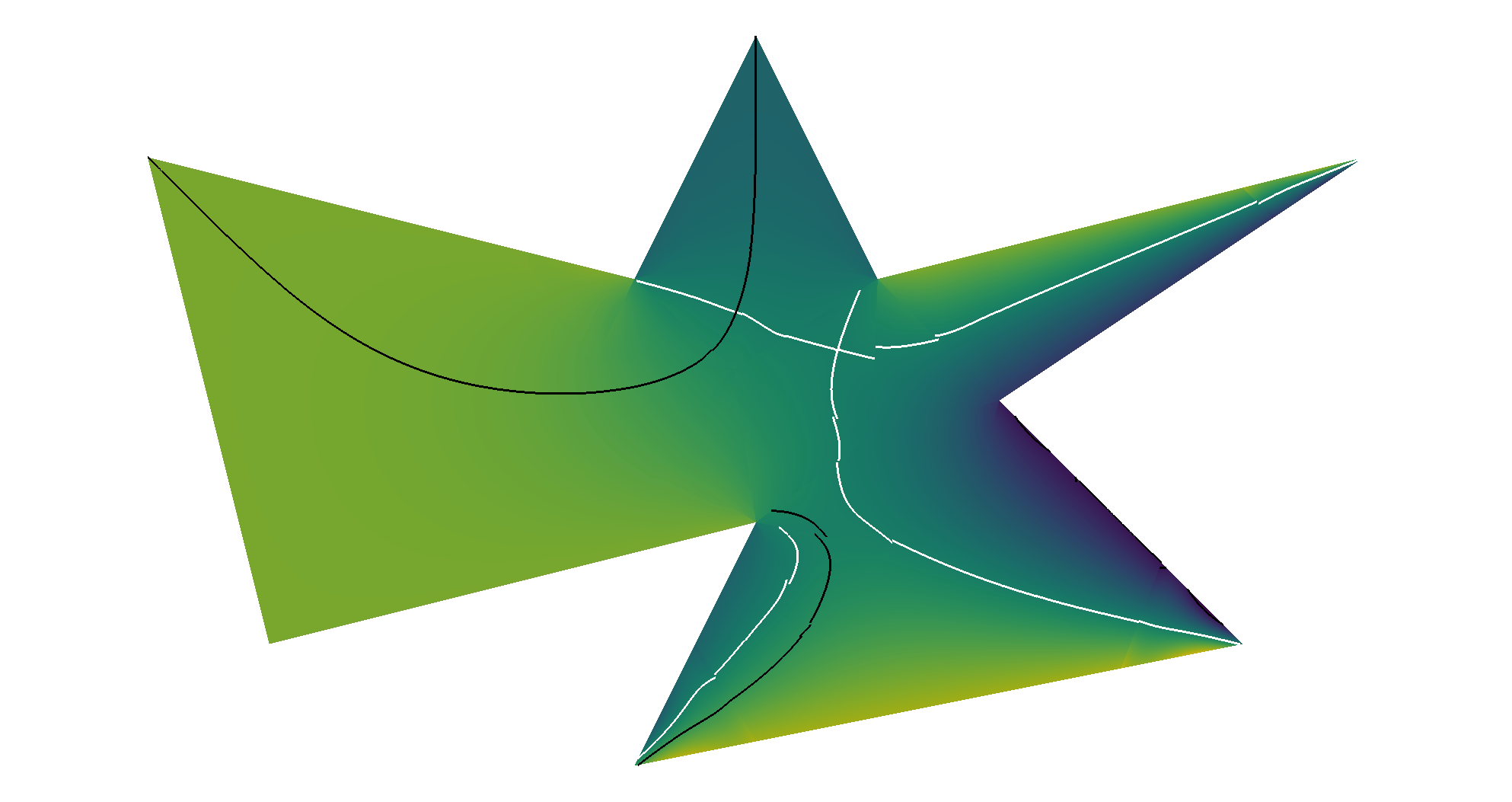}
  \includegraphics[width=0.49\textwidth]{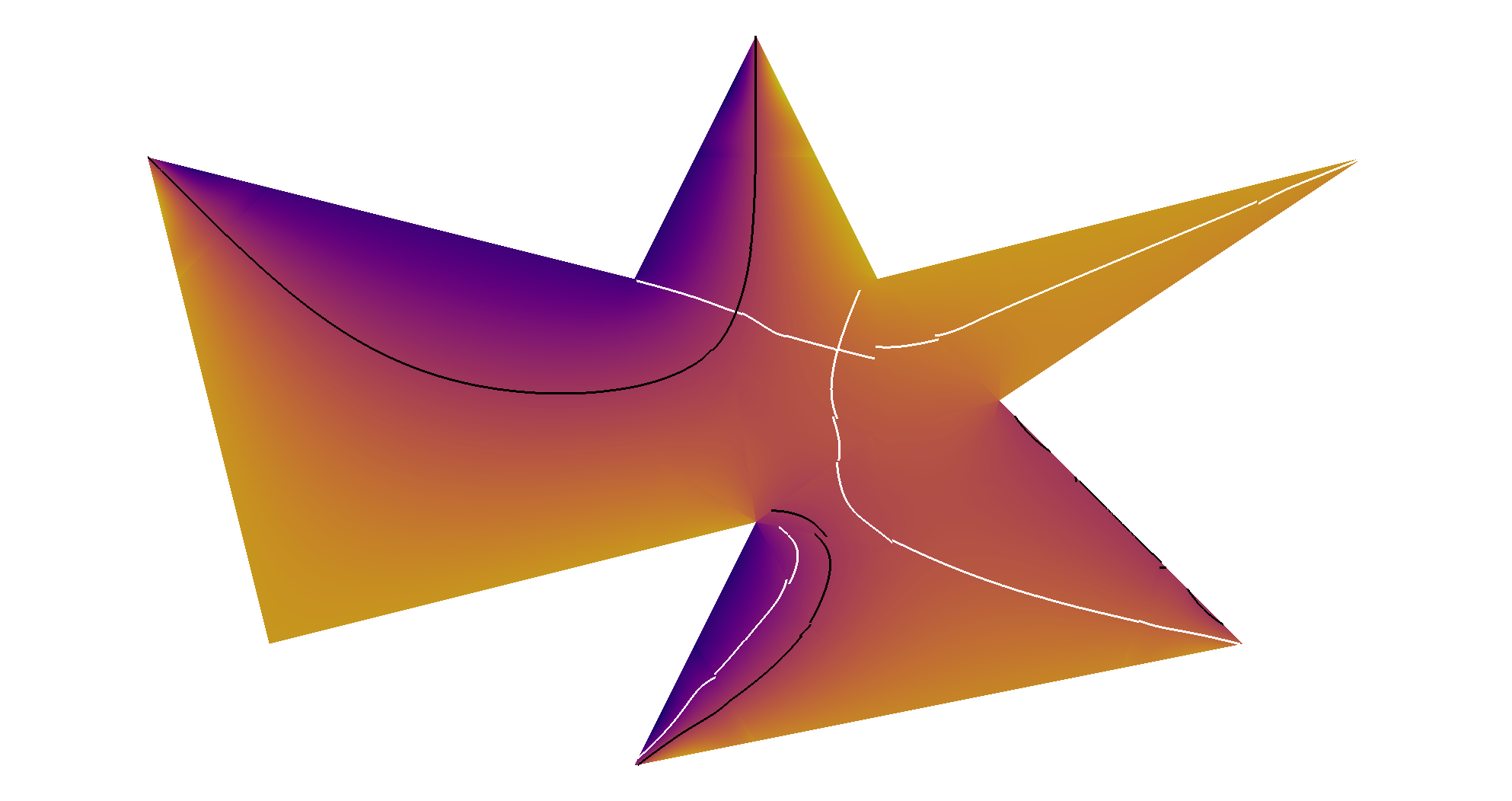}
  \includegraphics[width=\textwidth]{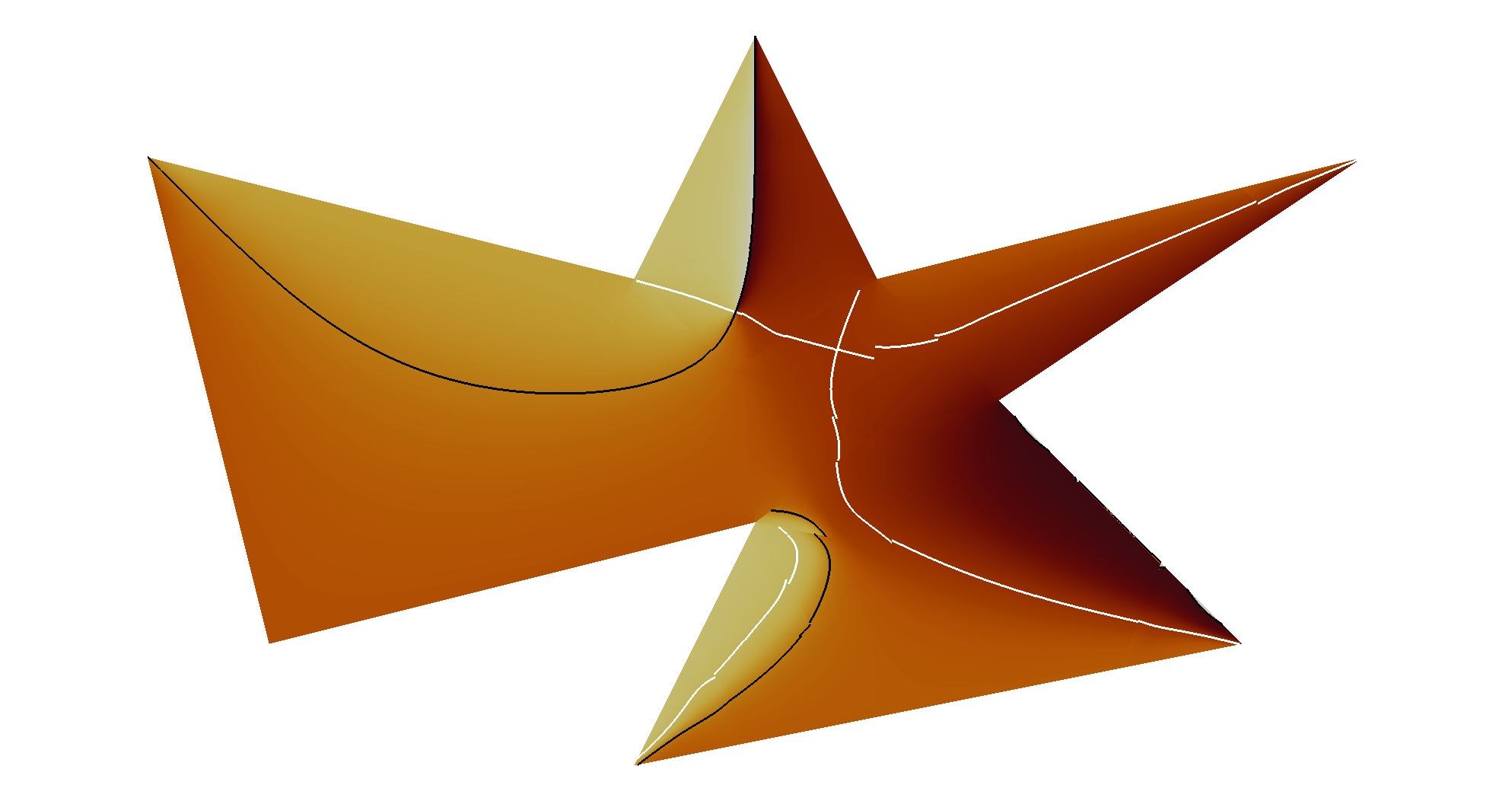}
  \caption{Solution \(u\) (top left) and \(v\) (top right) fields and computed \(\psi \) field (bottom) on a polygon with acute and obtuse corners using a DG discretization. Isocontours of \(u = 0\) and \(v = 0\) are shown in white and black respectively.}\label{fig:polygon-solution}
\end{figure}

As an example of a geometry that does not satisfy the smoothness conditions, we present the polygon geometry shown in Fig.~\ref{fig:polygon-solution}, which shows the solution and $\psi$ fields.

\begin{figure}[htbp]
  \centering
  \includegraphics[width=0.49\textwidth]{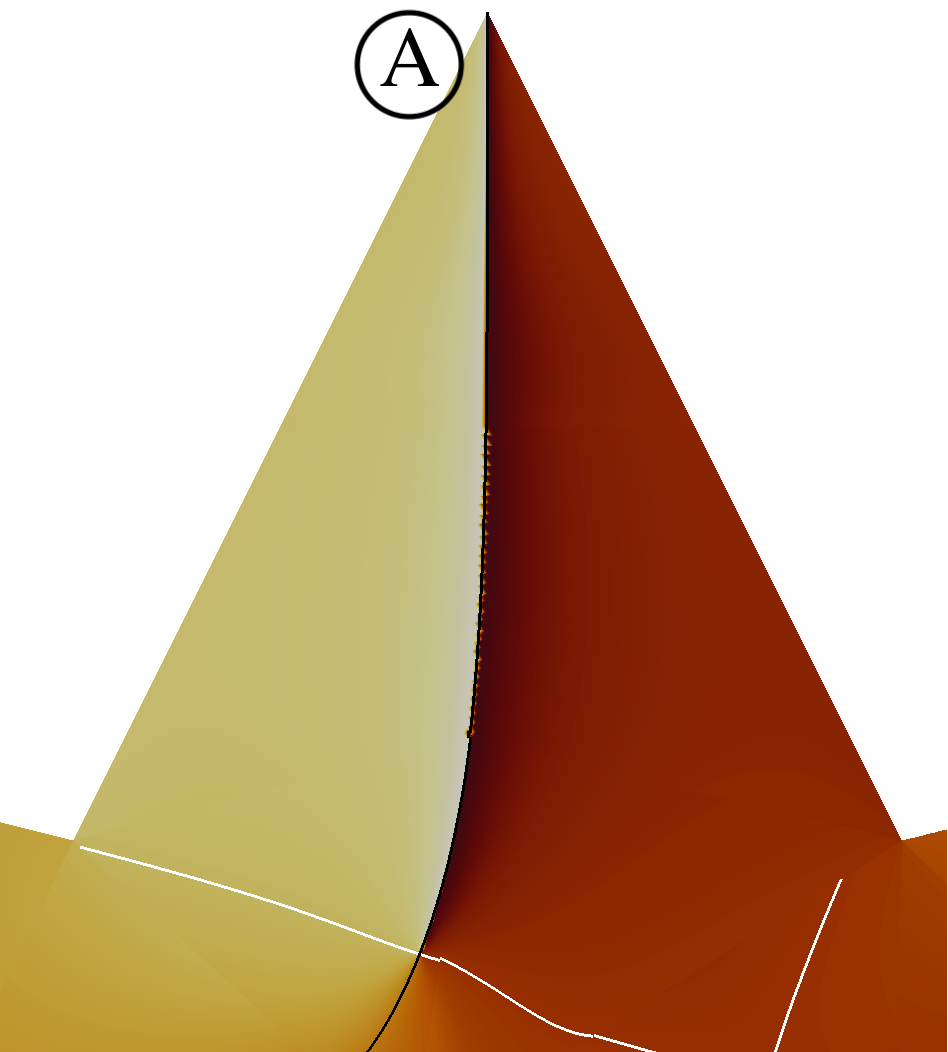}
  \includegraphics[width=0.49\textwidth]{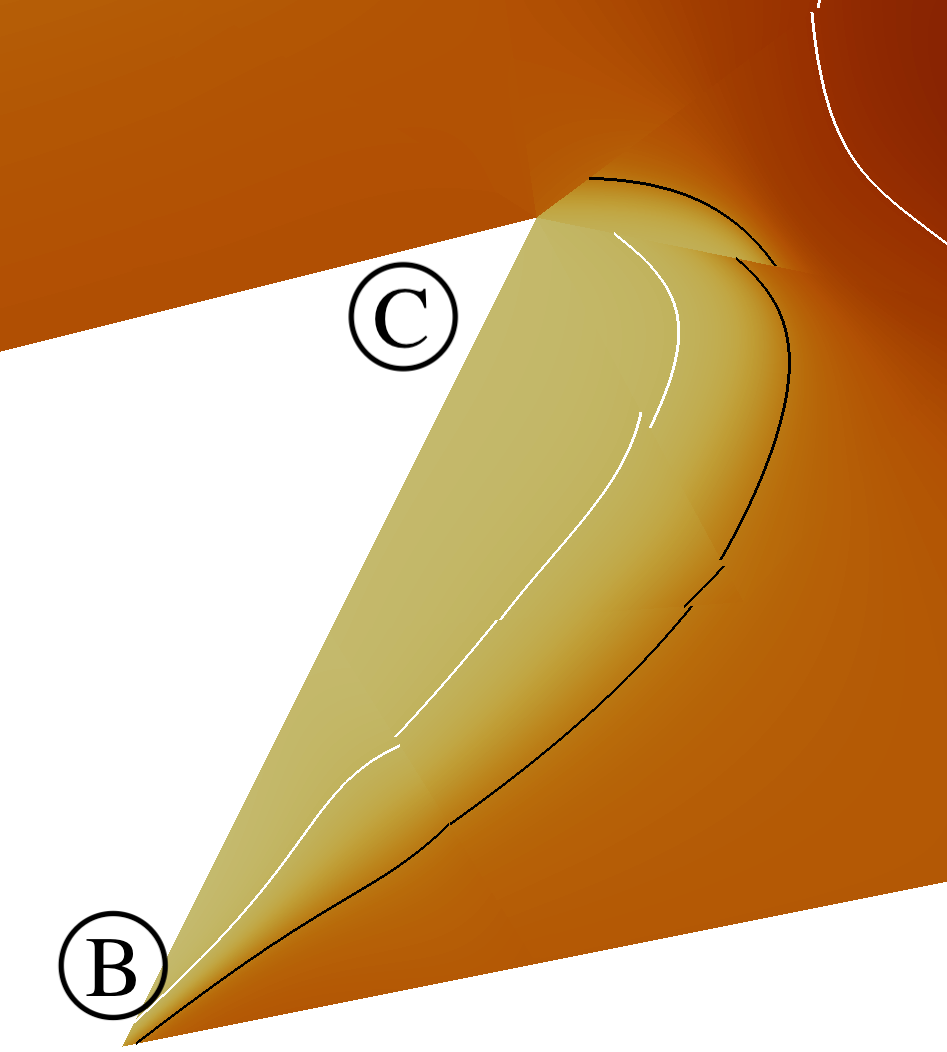}
  \caption{Closer views on the computed \(\psi \) field on the polygon geometry.}\label{fig:polygon-solution-close}
\end{figure}

We note two consequences of using a DG formulation to handle corners with angles that are not multiples of \(\pi/2\).
First, jump lines can originate from corners.
This fact is important, for it shows that one cannot \textit{a priori} and \textit{ad hoc} determine the valence of a boundary corner solely from its angle.
Fig.~\ref{fig:polygon-solution-close} (left) shows how the jump line originating from the top corner \circled{A} in the geometry ends in the sole {\color{black}critical point} of the domain.
Although we have observed that curvature produces {\color{black}critical points}, see Sec.~\ref{sec:math_critical}, this straight-sided example shows that a corner whose angle is not a multiple of \(\pi/2\) can also generate interior {\color{black}critical points}, again, see Fig.~\ref{fig:polygon-solution-close} (left).

Second, observe that discontinuous BCs are naturally enforced in the DG formulation.
No \textit{ad hoc} smoothing of the boundary is required, as shown in Fig.~\ref{fig:polygon-solution-close} (right), because the Dirichlet BCs can be enforced either through fluxes when one triangle element used to compute the guiding field shares both edges of the corner (see lower corner, \circled{B}), or exactly for all other corners (see upper corner, \circled{C}).

When a DG method is used to solve the Laplace problems for the guiding fields, the fields are discontinuous between elements.
This is especially true in domains where discontinuous BCs are present, as is clear in Fig.~\ref{fig:polygon-solution-close} (right).
Enough resolution (either through \textit{h}- or \textit{p}-refinement) needs to be provided so that the jumps at the element interfaces are small.
Too low resolution could result in zero isocontours located within solution jumps between elements.
The current implementation would not be able to detect such isocontours, nor could it consequently detect {\color{black}critical points}.
See Sec.~\ref{sec:critical} for more details.

Fortunately, {\color{black}critical points} should naturally be located in regions of smooth fields where jumps between elements are small.
Discontinuities between elements are large only near boundaries where discontinuous BCs appear.
The field on the inside of the domain is mostly smooth, due to the properties of the Laplace equation, resulting in very small jumps between elements.
If strong curvature is present, we observe that {\color{black}critical points} can appear near curved boundaries where BCs are continuous and the field is smooth.
If no (or little curvature) is present, {\color{black}critical points} are controlled by discontinuous BCs and are seen far away from the boundaries, on the inside of the domain, again where the field is smooth.
In all cases, {\color{black}critical points} are not expected to naturally occur in areas where large jumps between elements exist due to nearby discontinuous BCs.

Subsequently, the integration of streamlines is unaffected by the discontinuous nature of the discretization.
By using a multi step high order integrator, as explained later in Sec.~\ref{sec:streamlines}, streamline direction is not strongly affected by abrupt changes of direction between elements.
Streamlines remain smooth even when gradients in the guiding field are large, e.g.~when crossing a large jump between elements.

\subsection{{\color{black}Critical Point} Detection and Valence Evaluation}\label{sec:critical}

Once the solution to \eqref{eq:LaplaceEqns} is computed, zeros in the field \(\vec v\) are found and their valences are evaluated.
Unlike the usual cross field approach, we never, in fact, generate crosses.
As seen in Sec.~\ref{sec:MathematicalFormulation} the analysis can be performed on the \(\vec v\) field, and \(\psi \) is computed only when necessary, i.e.~for computing the valence of vertices and tracing streamlines via~\eqref{eq:math_corner} and~\eqref{eq:streamline-integration}.
All other operations can be accomplished by operating on \(\vec v\).

{\color{black}Critical} points inside the domain are located at \( \vec v = \vec 0 \).
To find those, we first search all elements and flag each that contains a \( \vec v = \vec 0 \) point at any quadrature point in the element.
If we find at least one value of each \( u \leq 0 \), \( u \geq 0 \), \( v \leq 0 \) and \( v \geq 0 \), a {\color{black}critical point} is either inside this or a nearby element.
In each flagged element, we then use Newton's method to approach the location of \( \vec v = \vec 0 \).
To do so efficiently, we perform the {\color{black}critical point} search in parametric space.
If a {\color{black}critical point} is found outside the reference element, we dismiss it and assume that it will be found through parametric search in a neighbor element.
In Fig.~\ref{fig:half-disc-solution}, we can clearly see the two symmetrically located {\color{black}critical points} where the black and white contours cross, consistent with {\color{black}the irregular nodes found by others}, e.g.~\cite{Nicolas-Kowalski:2012fu}, for this geometry.

After completing the search for all {\color{black}critical points}, a valence must be computed for each.
Referring to Sec.~\ref{sec:math_critical} and \eqref{eq:critical_valence}, we notice that the integral \(I_c\) can be evaluated solely by looking for the presence of jump lines.
In a counter-clockwise manner, a positive jump (i.e.~from \(-\frac{\pi}{4}\) to \(\frac{\pi}{4}\)) leads to a negative integral \( I_c = -1 \), and vice versa.
The sign of the jump can itself be determined solely by values of \(\vec v\), without the need to compute \(\psi \) or the cross field.

A jump must satisfy two conditions:
\(u\) must be negative and
\(v\) must change sign.
A positive increase of \(v\) indicates a positive jump, and vice versa.
To detect jumps, and therefore compute the valence of a {\color{black}critical point}, we step counter-clockwise over a sequence of uniformly distributed quadrature points located on a circle of small radius \(c\), centered on the {\color{black}critical point}.
Following this logic, we determine that the {\color{black}critical points} seen in Fig.~\ref{fig:half-disc-solution} both have a valence of three.

The valences of all boundary vertices are computed in a similar fashion, this time from expression~\eqref{eq:math_corner}.
The angle at the corner, determined from \(\theta_0\) and \(\theta_f\), is computed from the geometry.
To compute the integral \(I_c\) we compute \(\Delta\psi= \psi_f - \psi_0 \) from the values of \(\psi \) at the boundaries from the boundary conditions.
To \( \Delta\psi \), we add a jump contribution, if present, in the same way as for the computation of {\color{black}critical point} valences.
Because both corners in Fig.~\ref{fig:half-disc-solution} have an angle of \(\frac{\pi}{2}\), a valence of one can be trivially determined.

It is important to note that the boundary vertex operation does not require an \textit{a priori} determination of the corner valence based on angle only.
The corner valences can be determined through equation~\eqref{eq:math_corner} which does not only depend on the \textit{known} BCs but also on the existence of jump lines in the \textit{computed \(\psi \) field} in the neighborhood of the vertex, like that seen in Fig.~\ref{fig:polygon-solution-close}.

The boundary vertex valence gives the number of quads located at the associated corner.
Due to the definition, this valence could be zero.
Indeed, \( \Delta\theta \) and \( \Delta\psi \) could cancel out if the corner is sufficiently sharp.
Physically, this indicates the presence of a degenerate quadrilateral block (i.e.~a triangular block) where all streamlines would converge towards the degenerate corner.
This topology requires \textit{ad hoc} manipulation where the streamlines are used later to construct the quadrilateral subdivision of the domain.
This manipulation will be explained in the next section.

\subsection{Streamline Integration and Manipulation}\label{sec:streamlines}

Streamlines are traced throughout the domain after the valences of all {\color{black}irregular nodes} and vertices are determined.
The first step is to find the initial direction of each streamline.
Because \( \vec v = \vec 0 \) at {\color{black}critical points}, \(\psi \) is undefined.
Therefore, we choose to evaluate one of the streamline angles from \(\vec v\) at a small distance \(c\) of the {\color{black}critical point}.
This distance is currently set empirically based on the size of the elements.

The initial direction is refined iteratively from an initial guess using Algorithm~\ref{alg:InitialStreamlineDirection}, which is inspired by bisection.
For the first streamline, we can take any initial guess for the direction.
For subsequent streamlines, we take an initial guess at angles multiples of \( \frac{2\pi}{\mathcal V} \).
The size of the tolerance, $\epsilon$ is not critical; we use $10^{-9}$.

\begin{algorithm}[htbp]
  \caption{Determination of the initial direction of a streamline at an {\color{black}irregular node} $\vec p_{0}$.}\label{alg:InitialStreamlineDirection}

  \KwData{Direction \( \alpha_0 \) (initial guess), {\color{black}Irregular node} location \(\vec p_0\)}
  \KwResult{Direction \( \alpha_f \) (converged value)}

  Initialize \( \alpha \) from \( \alpha_0 \) \;
  \While{\(\left| \Delta\alpha\right| > \epsilon \) }{
    Compute point \(\vec p_1\) at distance \(c\) and direction \( \alpha \) of \(\vec p_0\) \;
    Interpolate \(\vec v\) and compute \( \psi \) at \( \vec p_1 \) \;
    Find \( \psi' = \psi + k\frac{\pi}{2} \) where \( k = 0,1,2,3 \) such that \( \Delta\alpha = \left| \alpha - \psi' \right| \) is minimized \;
    Update \( \alpha = \psi' \) \;
  }
  Return \( \alpha_f = \alpha \) \;
\end{algorithm}

The angle search in Alg.~\ref{alg:InitialStreamlineDirection} can also be applied to boundary vertices, including those where the boundary conditions might be discontinuous and \(\psi \) is ambiguous.
An initial guess for \textit{all} streamlines may be obtained at angles multiples of \( \frac{\Delta\theta}{\mathcal V} \).

After an initial direction is obtained, streamlines are \emph{synchronously} advanced throughout the domain.
This part of the procedure is the most computationally expensive.
Each new streamline point requires the search for the element that contains it.
Next, the inverse map \( \vec \xi = \mathcal{X}^{-1}_e( \vec x )\) is used to transform to parametric coordinates.
Finally the high order interpolation of \(\vec v\) is computed via eq.~\eqref{eq:SolutionInterpolant}.

A multi stage integration like Runge-Kutta, where intermediate evaluations are used, would prove too expensive to integrate the streamlines.
Instead, we use a \(4^{th}\) order multi step Adams-Bashforth integrator.
The current implementation reduces the order of the algorithm when too few steps are available.
While advancing the streamlines, it is important to validate the direction \( \alpha \) obtained at the latest point.
When a streamline crosses a jump line, \( \psi \) abruptly rotates \( \pm\pi/2 \) and \( \alpha \) must be adjusted accordingly.

As streamlines are advanced throughout the domain, they may meet.
To anticipate two streamlines meeting, the front points of each are compared at each step.
If the distance between the front points of two streamlines is less than the step size and they are advancing in opposite directions, the streamlines are assumed to meet.
When the absolute value of the difference in the latest $\alpha$ values equals $\pi$ when rounded to the nearest $\pi/2$, we consider the streamlines to be advancing in opposite directions.
Each of two meeting streamlines are then advanced to the starting point of the other, keeping the same number of points in each.

A merged streamline is created by weight-averaging a pair of streamline points.
A simple linear weight function could be used, but that will change the angles computed with Alg.~\ref{alg:InitialStreamlineDirection} at the starting points, both {\color{black}irregular nodes} and boundary vertices.
This not desirable as the PDE problem has already optimized these angles to be as equally distributed as possible.
Instead we use trigonometric weight functions of the form \( W_0(x) = \cos^2(x) \) and \( W_1(x) = \sin^2(x) \), where $x\in[0,1]$ is the parametrization along the streamline.
Because \( \frac{dW_0}{dx} \) and \( \frac{dW_1}{dx} \) are zero at \( x=0 \) and \( x=1 \), angles are preserved at the initial and final points of the merged streamline.

{\color{black}
\subsubsection{Aggressive Merging}

What has been described above constitutes expected or \emph{normal} merging.
But the user may also want to reduce the number of separatices and thereby simplify a valid block decomposition.
One approach is to perform \emph{aggressive} merging of the separatrices.
For aggressive merging, a larger distance threshold is carefully chosen to make streamlines that would normally just miss each other merge anyway.
Trigonometric weight functions are especially useful when performing aggressive merging of two reasonably distant streamlines.
An example of aggressive merging with trigonometric weight functions is presented in Fig.~\ref{fig:nautilus-streamlines} where the left figure corresponds to normal merging and the right one to aggressive merging.
In this case, a distance threshold of 5 times the step size was used for aggressive merging.
This example shows that the separatrix graph is greatly simplified and that the curving of some streamlines helps preserve the original direction at each end.
Note, importantly, that irregular nodes obtained from the critical points are not modified \textit{a posteriori};
only separatrices are.

\begin{figure}[htbp]
  \centering
  \includegraphics[trim={30cm 0 30cm 0},clip,width=0.49\textwidth]{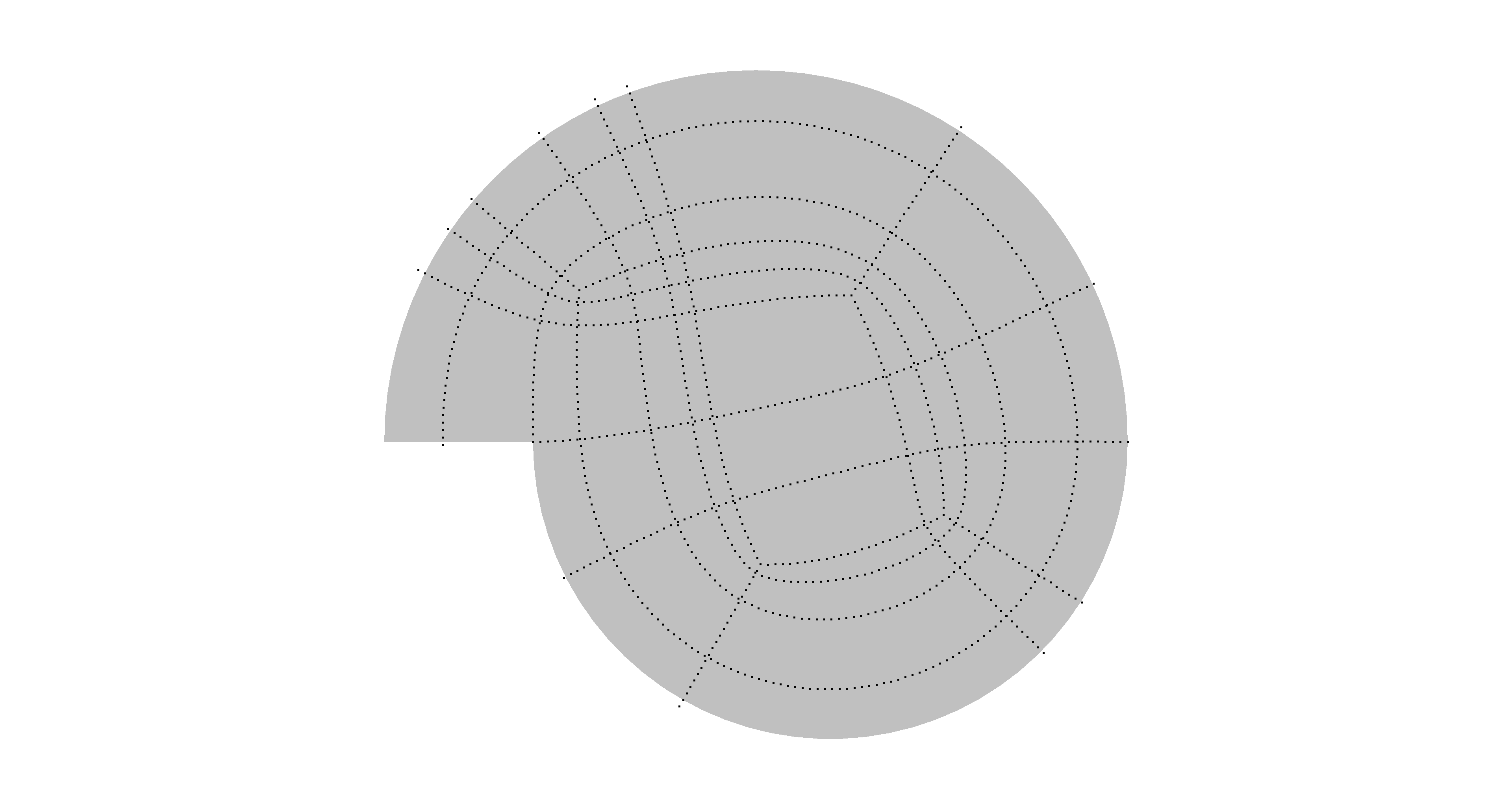}
  \includegraphics[trim={30cm 0 30cm 0},clip,width=0.49\textwidth]{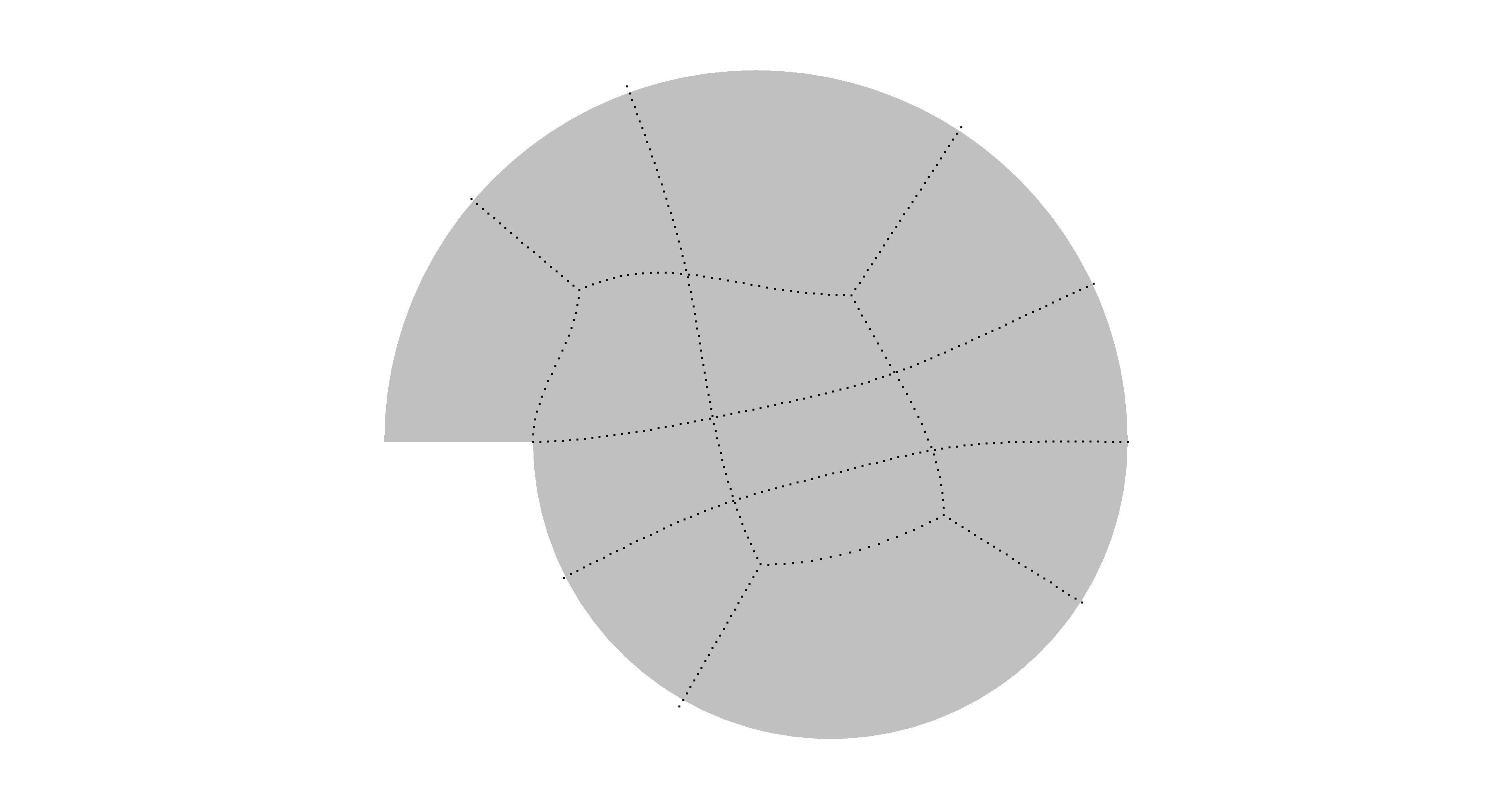}
  \caption{Example of normal (left) and aggressive (right) merging. A distance threshold of 5 times the step size was used for aggressive merging.}\label{fig:nautilus-streamlines}
\end{figure}
}

\subsection{Quad Mesh Generation}\label{sec:quad-meshing}

\emph{NekMesh} relies on the OpenCASCADE platform~\cite{OpenCascadeSAS2018} as its CAD engine for mesh generation and for the projection to high order boundary representations.
OpenCASCADE also includes tools for CAD manipulation.
We make use of these capabilities for spline and wire creation and for face splitting.

\begin{figure}[htbp]
  \centering
  \includegraphics[width=0.49\textwidth]{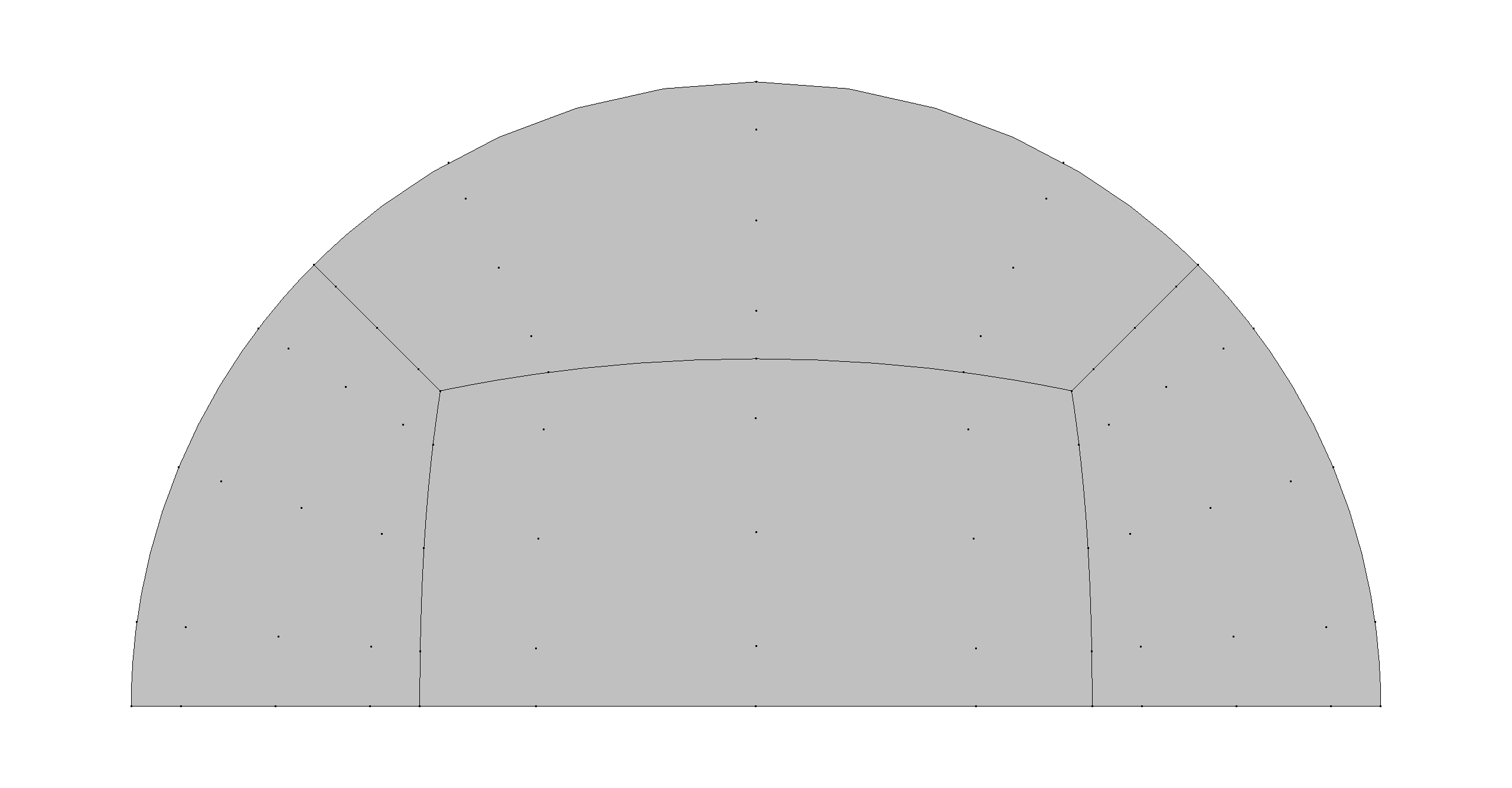}
  \includegraphics[width=0.49\textwidth]{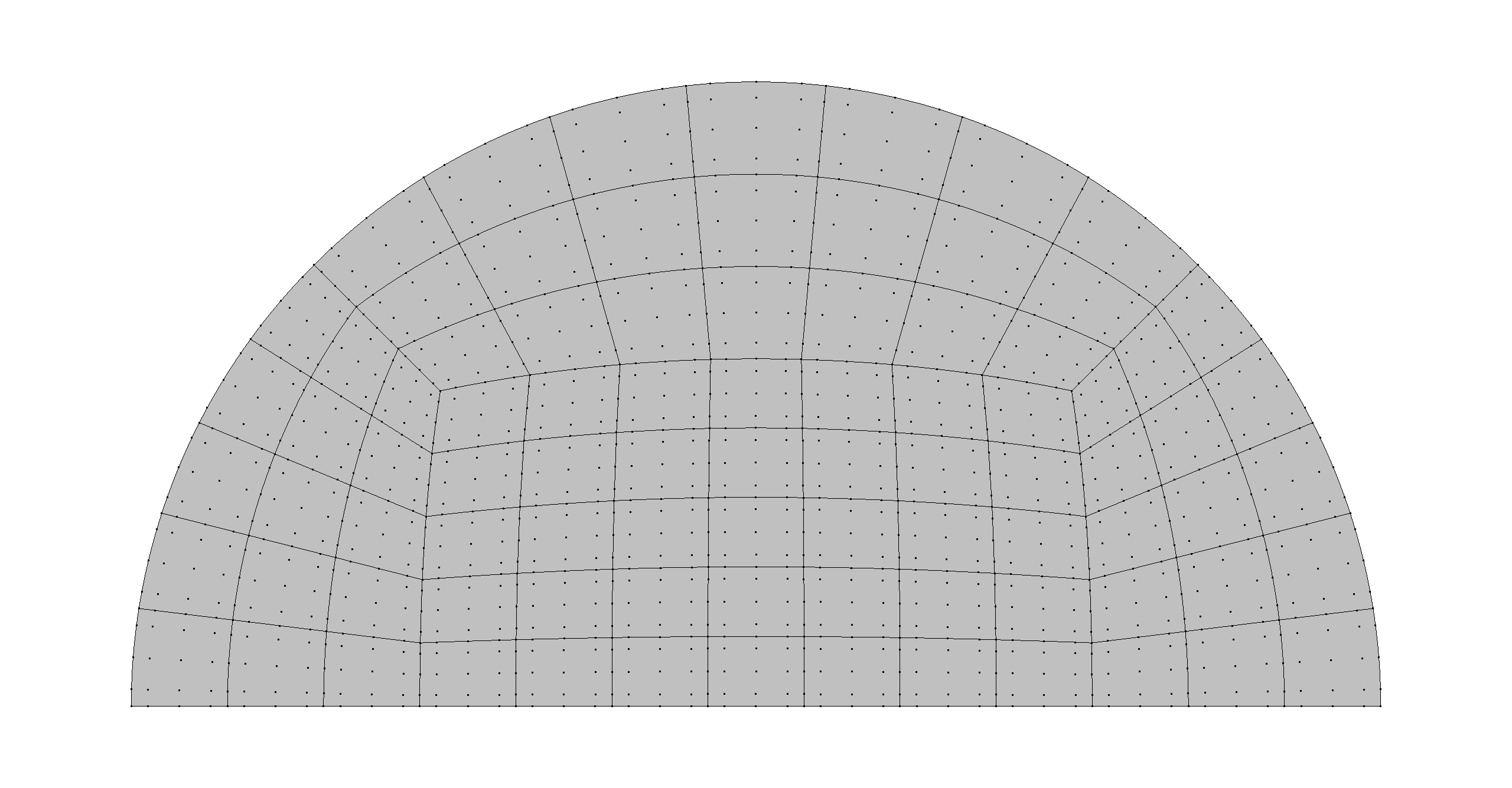}
  \caption{Coarse (left) and subdivided (right) quadrilateral meshes obtained on the half disc geometry based on detected {\color{black}irregular nodes} and traced and merged {\color{black}separatrices}.}\label{fig:half-disc-meshes}
\end{figure}

The {\color{black}separatrices} computed from the guiding field are transformed into OpenCASCADE \textit{edges} represented by interpolating splines.
These \textit{edges} can be joined into a set of \textit{wires}, with each \textit{wire} consisting of topologically connected \textit{edges}.
The original CAD file is loaded again, consisting of a single two dimensional \textit{face}, representing the domain $\Omega$.
This \textit{face} is then iteratively split by each {wire} and a set of quadrilateral \textit{faces} is obtained.
These faces are topologically connected, meaning that the future mesh will be conforming.
At this point, each {face} can be trivially meshed with a single quadrilateral and projected to high order.
Fig.~\ref{fig:half-disc-meshes} (left) illustrates the quadrilateral decomposition and coarse quadrilateral mesh obtained for the half circle after computing the {\color{black}location of irregular nodes} and tracing the {\color{black}separatrices}.

\begin{figure}[htbp]
  \centering
  \includegraphics[width=0.75\textwidth]{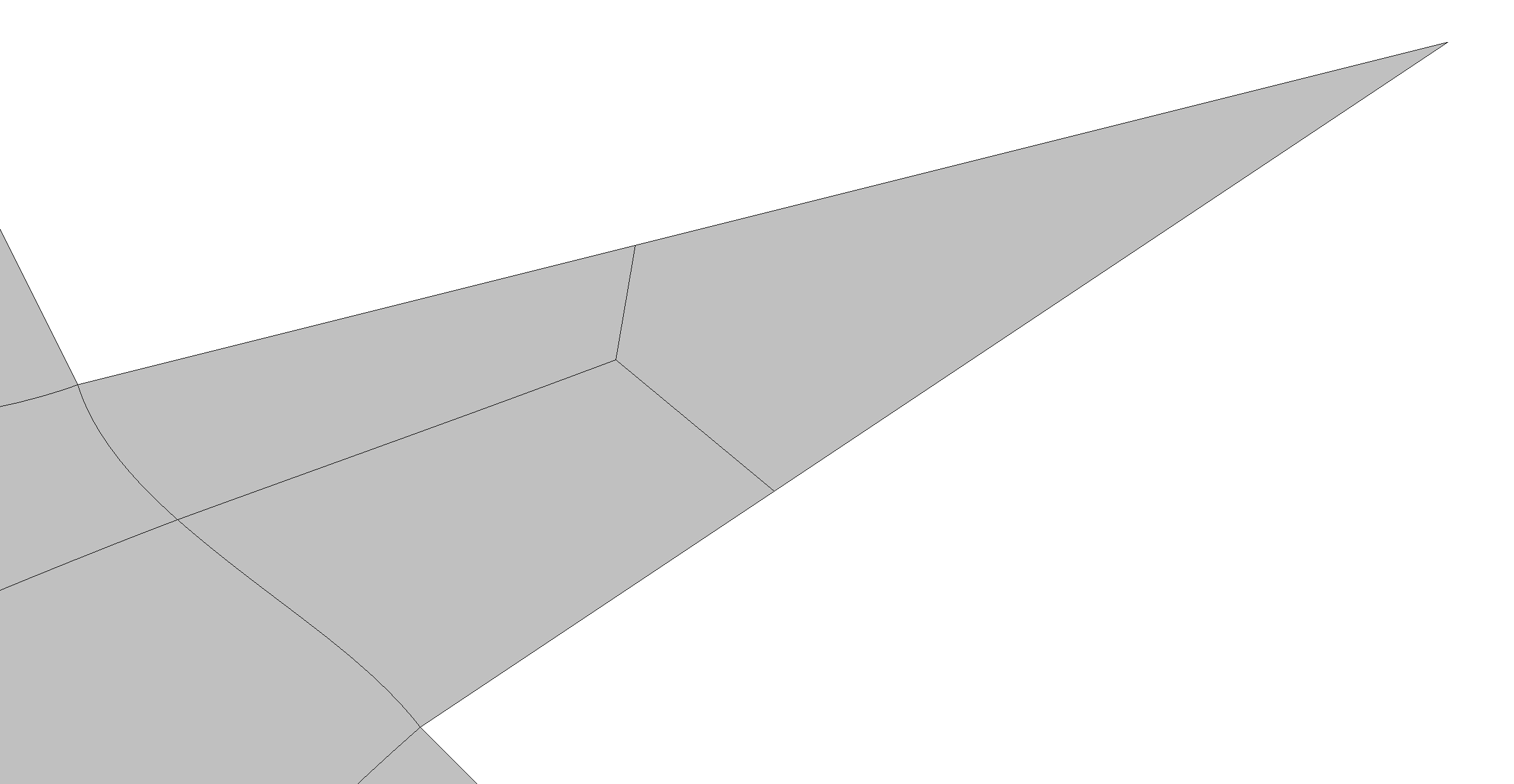}
  \caption{Detail of a domain showing a midpoint division of a degenerate quadrilateral caused by a sharp angled corner. The full domain is shown in Fig.~\ref{fig:polygon-mesh}.}\label{fig:midpoint-division}
\end{figure}

We use the midpoint division approach~\cite{Li1995} to split a triangular block that forms when a corner valence is zero.
We insert an artificial 3-valence {\color{black}irregular node} which we connect to each of the three sides of the triangle.
One {\color{black}separatrix} can then be physically integrated away from the degenerate corner and throughout the domain.
The other two branches of the 3-valence {\color{black}irregular node} are defined as straight lines at \( \pm\frac{2\pi}{3} \) from the physical {\color{black}separatrix}.
Fig.~\ref{fig:midpoint-division} shows what the midpoint division of a degenerate quadrilateral might look like.

If desired, we make use of an isoparametric splitting approach described in~\cite{Moxey2015} to further subdivide the quadrilateral subdomains into smaller elements.
This technique, initially developed for boundary layer mesh division, guarantees that subdivided elements have the same quality as the macro element from which they are created.
In the current context, it also ensures that angles are preserved.
Each row of quads can be individually split using the isoparametric approach.
The right-hand side of Fig.~\ref{fig:half-disc-meshes} shows what the mesh on the left looks like after a conformal splitting of each row of quads.

\section{Examples}\label{sec:examples}

We present four examples of meshes generated for geometries taken from, or inspired by, the literature.
The meshes shown in Fig.~\ref{fig:half-disc-meshes} are typical of those already generated by low order approaches~\cite{Viertelabs-1708-02316,Nicolas-Kowalski:2012fu,Viertel:2017rt}.
We use the additional examples to illustrate the application of the high resolution approach on multiply connected geometry, one for which low order approaches have been shown to have difficulty, a polygon with both acute and obtuse angles, and a NACA 0012 airfoil.

The geometry for Example \rm{I} shown in Fig.~\ref{fig:iti-spurious} is a rectangle with two quarter circle holes.
This geometry has also been meshed, for example, in reference~\cite{Nicolas-Kowalski:2012fu}.
{\color{black}
Some low order implementations have produced spurious asymmetric irregular nodes, see Fig.~\ref{fig:iti-spurious}~\cite{iti}, for example.
The high order guiding field approach, however, produces the minimum number of {\color{black}irregular nodes} (two 3-valence {\color{black}irregular nodes}) while preserving the symmetry of the geometry without the need to coalesce spurious {\color{black}irregular nodes} that can result from low order detection algorithms.
Fig.~\ref{fig:iti-mesh} shows the coarse quadrilateral mesh obtained on this geometry using the guiding field approach.

\begin{figure}[htbp]
  \centering
  \includegraphics[width=0.75\textwidth]{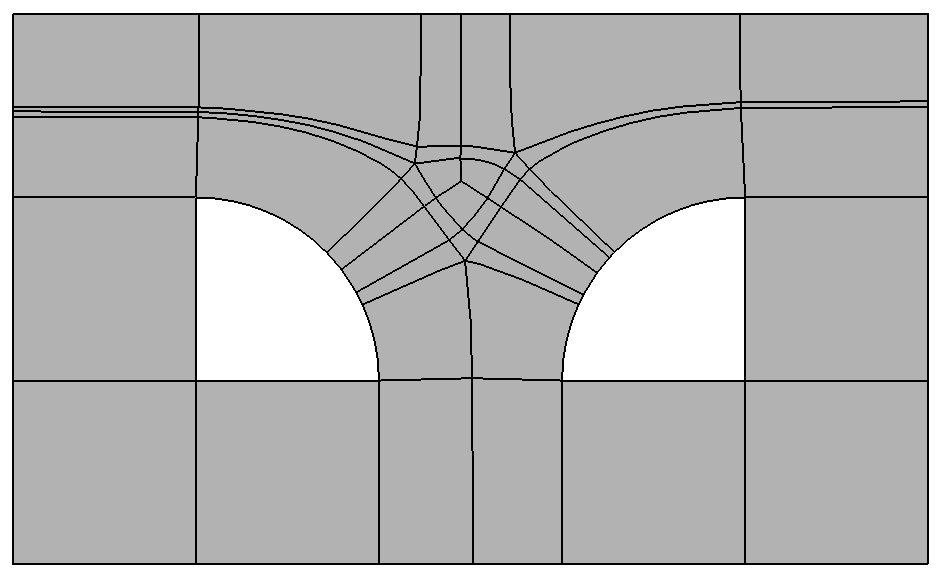}
  \caption{Block decomposition of Geometry \rm{I} with spurious asymmetric irregular nodes, generated with a low order cross field implementation~\cite{iti}.}\label{fig:iti-spurious}
\end{figure}

\begin{figure}[htbp]
  \centering
  \includegraphics[width=\textwidth]{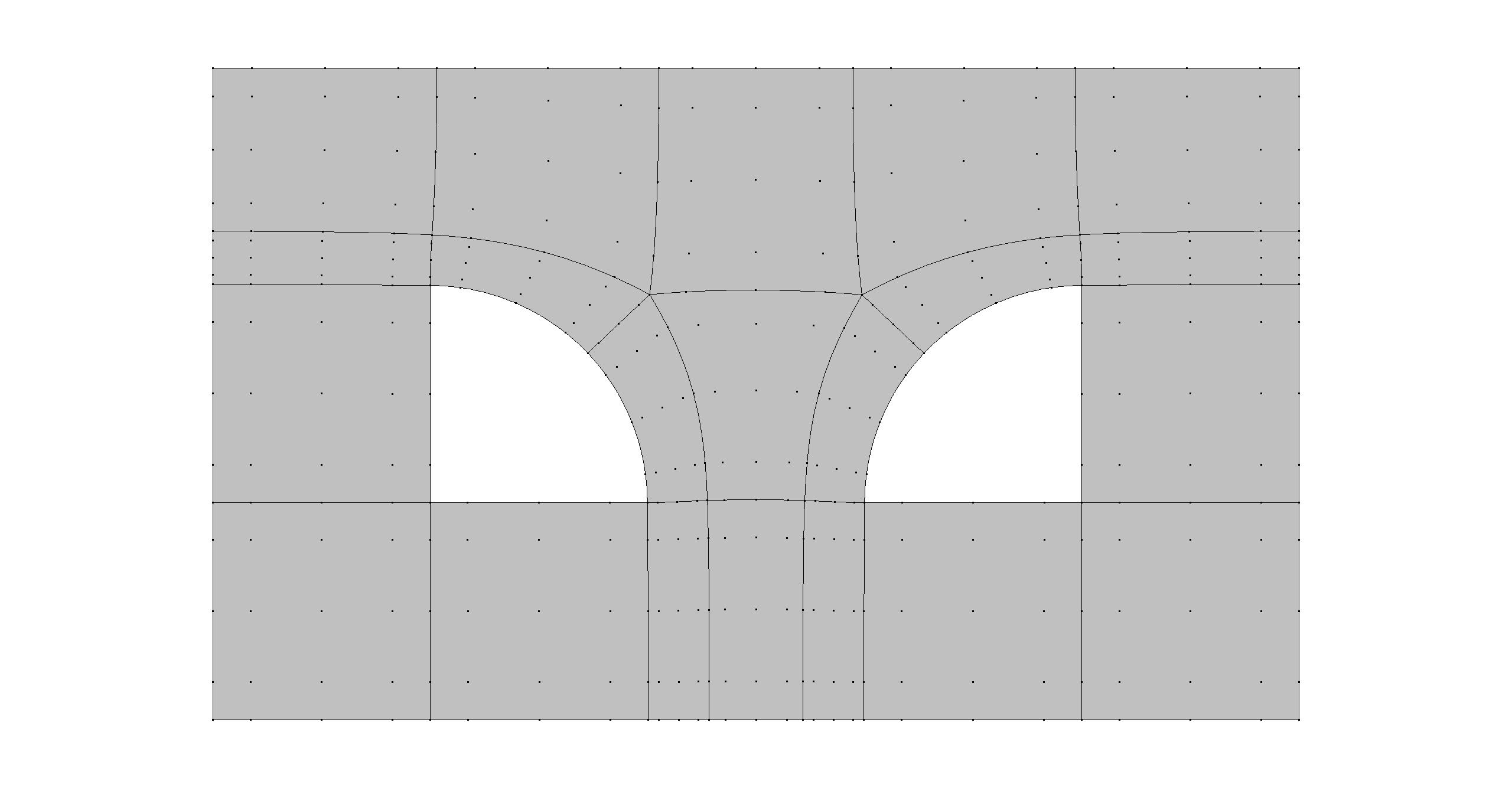}
  \caption{Coarse quadrilateral mesh obtained on Geometry \rm{I} using the guiding field approach.}\label{fig:iti-mesh}
\end{figure}
}

\begin{figure}[htbp]
  \centering
  \includegraphics[width=\textwidth]{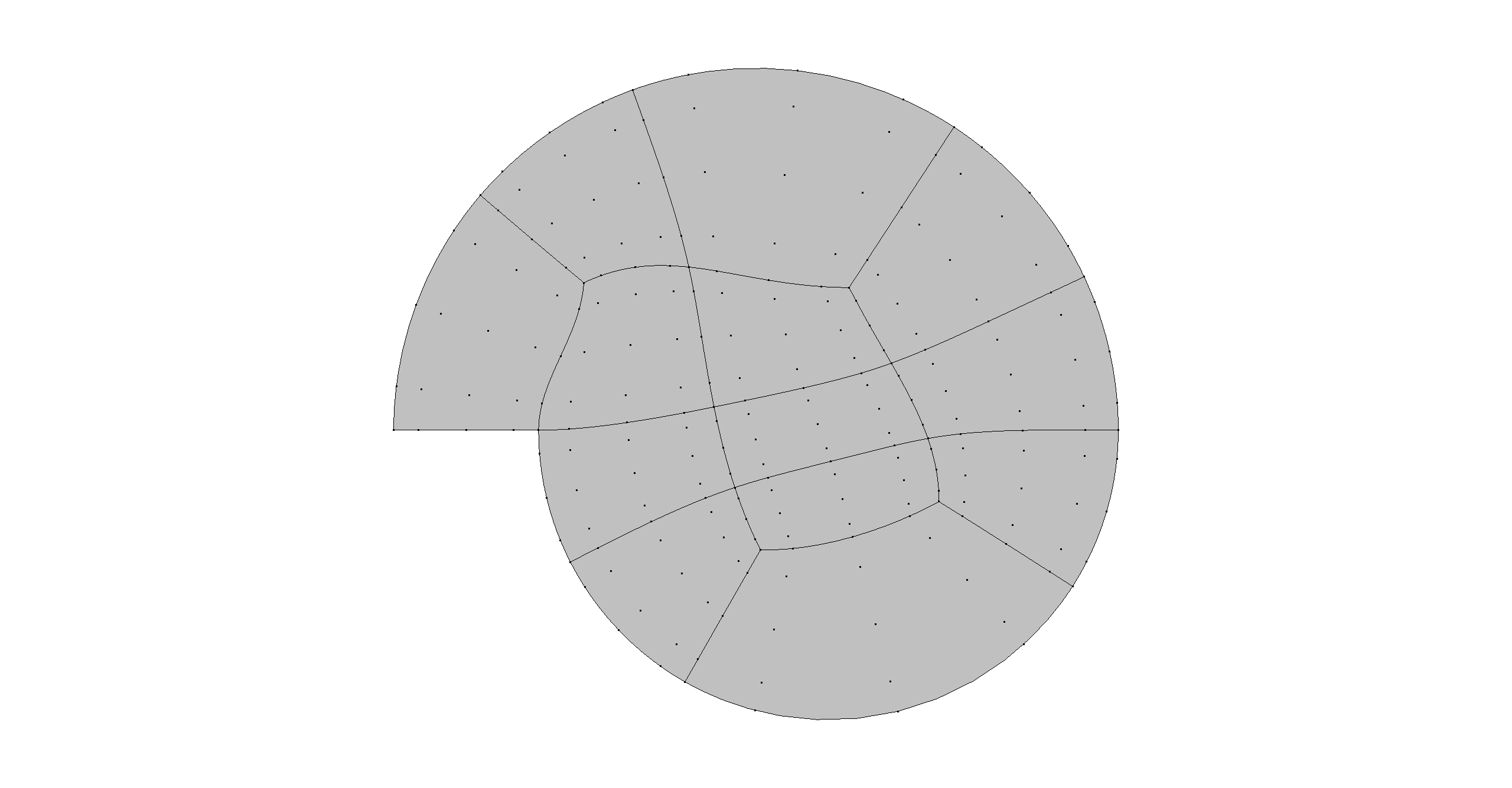}
  \caption{Coarse quadrilateral mesh obtained on Geometry \rm{II}, the nautilus.}\label{fig:nautilus-mesh}
\end{figure}

{\color{black}Example \rm{II}, is the nautilus.
The domain and guiding fields shown in Fig.~\ref{fig:nautilus-solution} give the decomposition shown in Fig.~\ref{fig:nautilus-mesh} when agressive merging is applied.
This domain is one that low order methods generate spurious spiral separatices that end in a limit cycle, therefore creating an invalid block decomposition~\cite{Viertelabs-1708-02316,Viertel:2017rt}.
To date, no approach using cross fields has been presented that generates a valid separatrix graph.
The use of high order integration of an accurate guiding field allows \emph{NekMesh} to generate {\color{black}separatrices} that do not form a spiral.
In fact, the {\color{black}irregular nodes} found by \emph{NekMesh} on this geometry are identical to those one would obtain for a simple disc.
Fig.~\ref{fig:nautilus-mesh} shows the coarse quadrilateral mesh obtained on this nautilus geometry.
This decomposition was generated using aggressive merging with a distance threshold of value 5 times the step size.
If we were to use normal merging, the block decomposition would have looked more complicated, as shown in Fig.~\ref{fig:nautilus-streamlines}~(left), reducing the flexibility for users to split the mesh to their liking.

We note that the nautilus has a similar (mirrored) guiding field to that obtained on the half disc in Fig.~\ref{fig:half-disc-solution}.
The nautilus itself consists of two half discs of different size attached by their chord.
The irregular node pattern (four 3-valence nodes) is therefore expected.
We also note one major difference with the irregular nodes obtained with previous techniques, e.g.~references~\cite{Viertelabs-1708-02316,Viertel:2017rt}:
The upper left irregular node is located further outside than in previous works.
This seems to lead to streamlines successfully escaping a limit cycle.

\begin{figure}[htbp]
  \centering
  \includegraphics[trim={20cm 0 20cm 0},clip,width=0.49\textwidth]{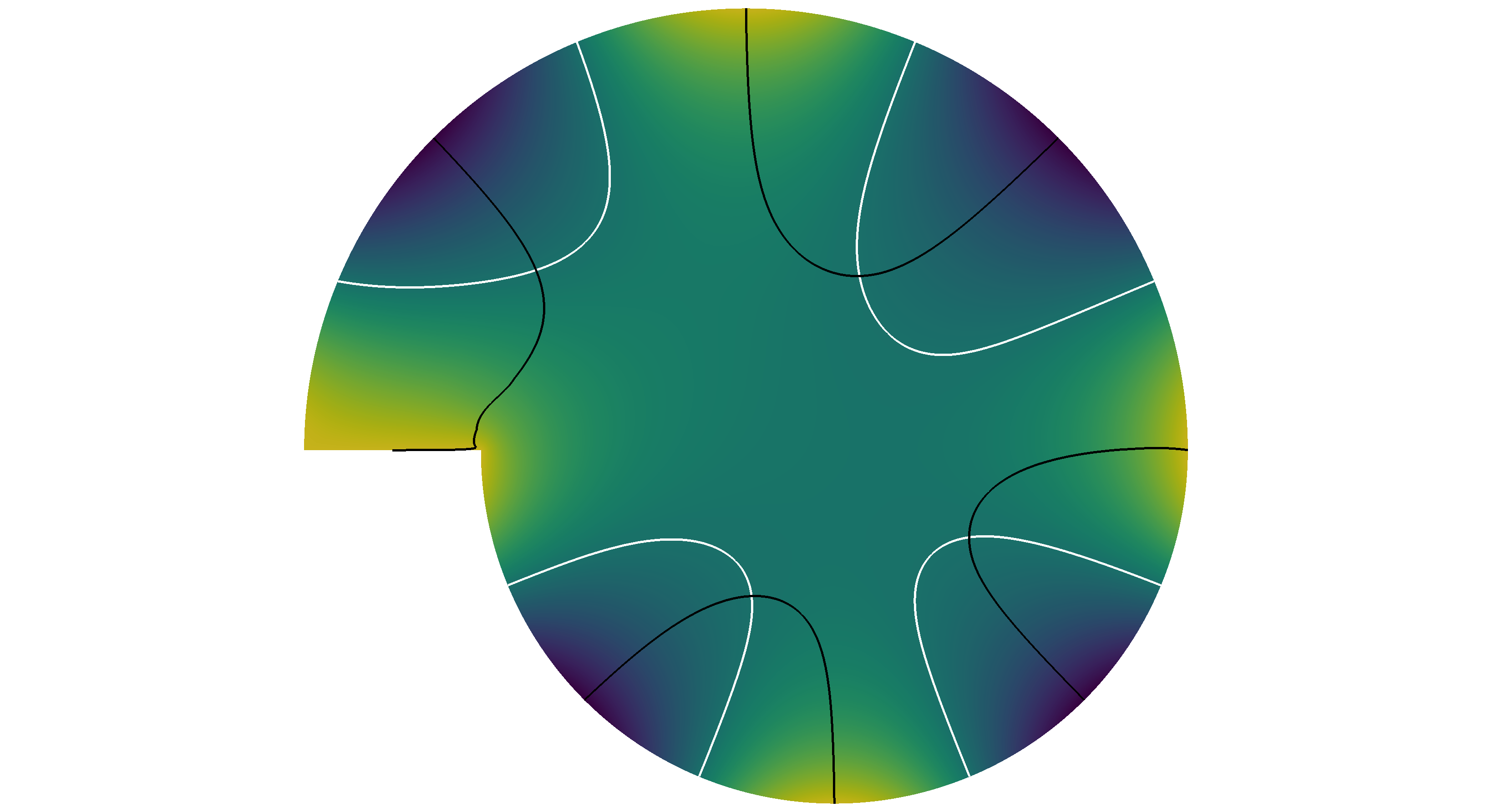}
  \includegraphics[trim={20cm 0 20cm 0},clip,width=0.49\textwidth]{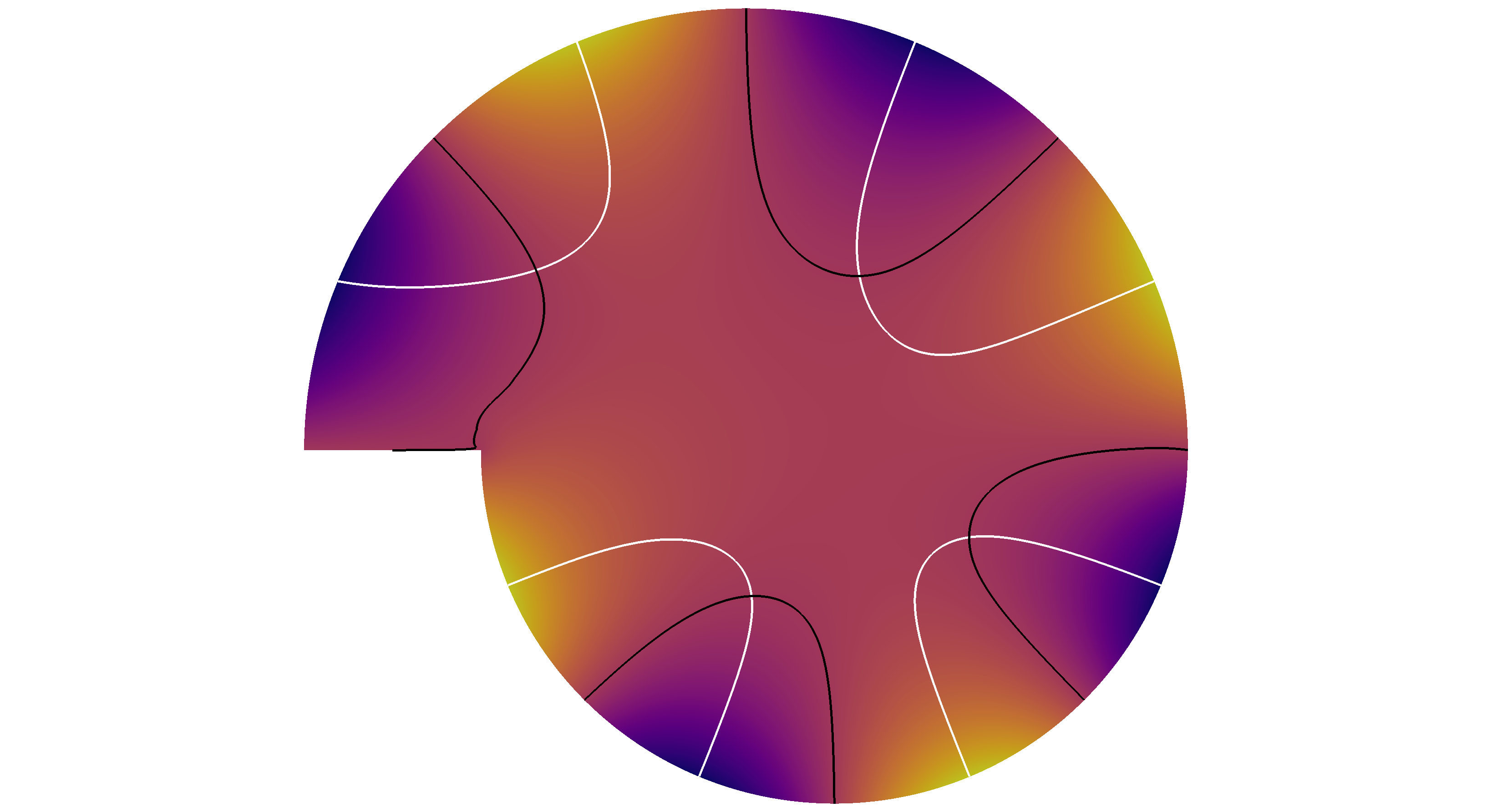}
  \includegraphics[trim={20cm 0 20cm 0},clip,width=0.75\textwidth]{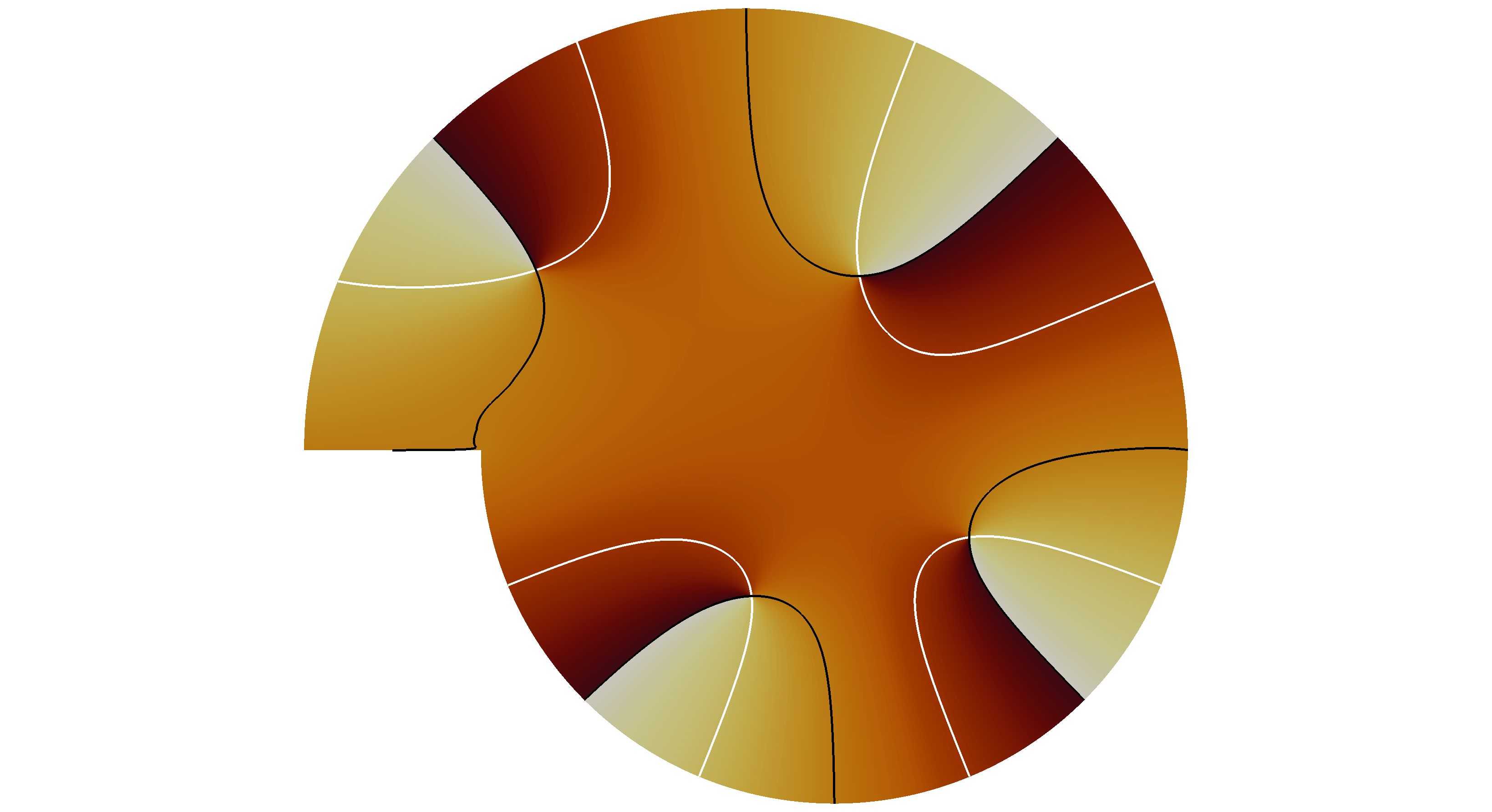}
  \caption{Solution \(u\) (top left) and \(v\) (top right) fields and computed \(\psi \) field (bottom) on the nautilus geometry. Isocontours of \(u = 0\) and \(v = 0\) shown in white and black respectively.}\label{fig:nautilus-solution}
\end{figure}
}

Example \rm{III} consists of a polygon whose corners include acute and obtuse angles.
We use this example to demonstrate the use of the DG discretization for the solution of the guiding field.
The rest of the meshing procedure remains the same and is in fact unaffected by the type of discretization.
Fig.~\ref{fig:polygon-solution} has already shown the solution \(\vec v\) and computed \(\psi \) field for this geometry, with
\(u = 0\) and \(v = 0\) isocontours shown in white and black, respectively.
As mentioned earlier in Sec.~\ref{sec:streamlines}, the only consideration when using the DG formulation is that the approximate solutions are accurate enough to not affect the streamline integration.

\begin{figure}[htbp]
  \centering
  \includegraphics[width=\textwidth]{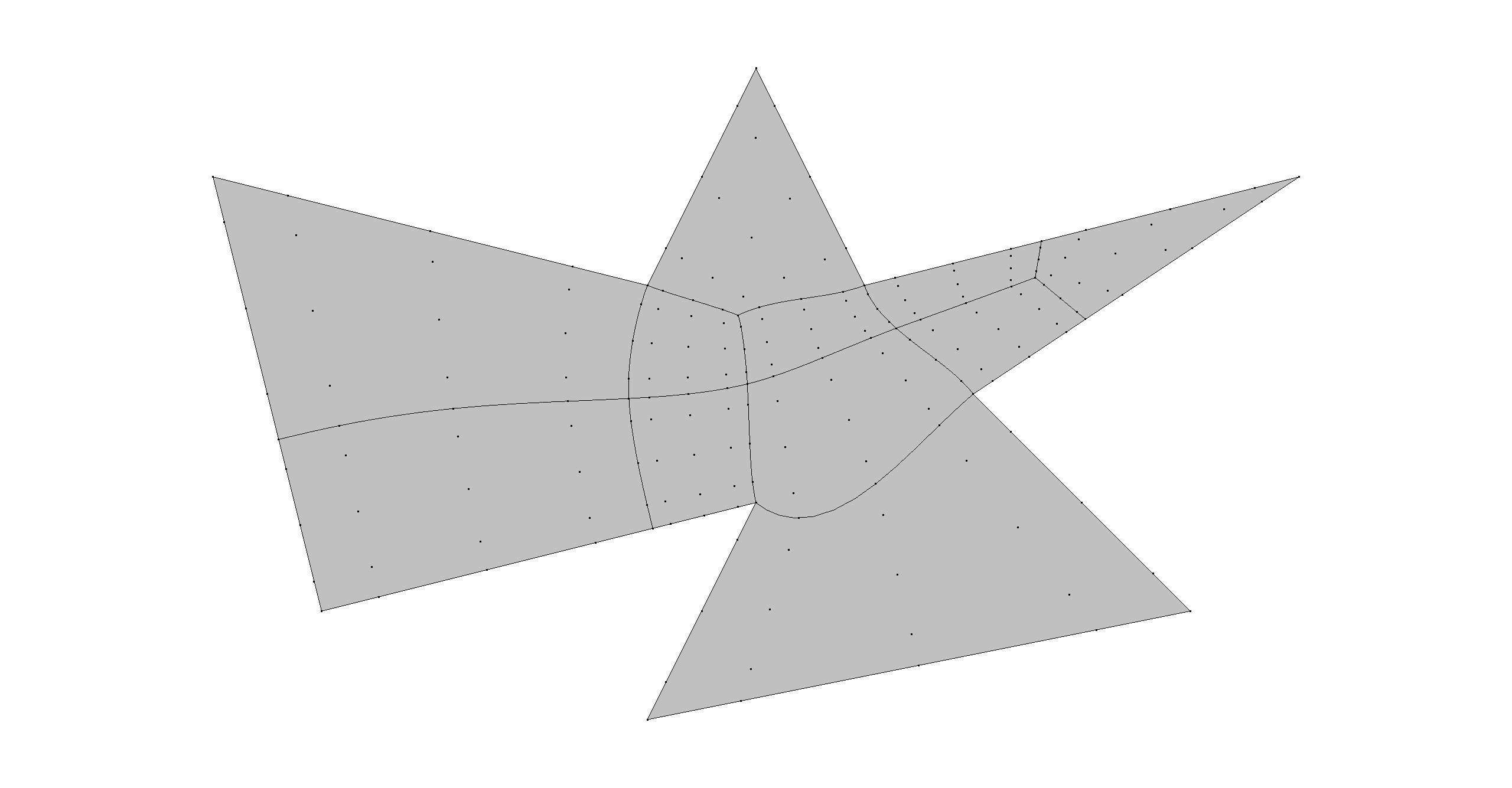}
  \caption{Coarse quadrilateral mesh obtained on Geometry \rm{III}, a polygon.}\label{fig:polygon-mesh}
\end{figure}

The Example \rm{III} geometry features a sharp corner (far right, Fig.~\ref{fig:polygon-mesh}) whose valence is evaluated as zero.
That creates a degenerate quadrilateral that must be split into three valid quads as explained in Section~\ref{sec:streamlines}.
Fig.~\ref{fig:polygon-mesh} shows the coarse quadrilateral mesh obtained on the polygon geometry.
It contains only two irregular nodes, one created \textit{ad hoc} for the degenerate corner and the other detected as a {\color{black}critical} point in the guiding field.

\begin{figure}[htbp]
  \centering
  \includegraphics[width=\textwidth]{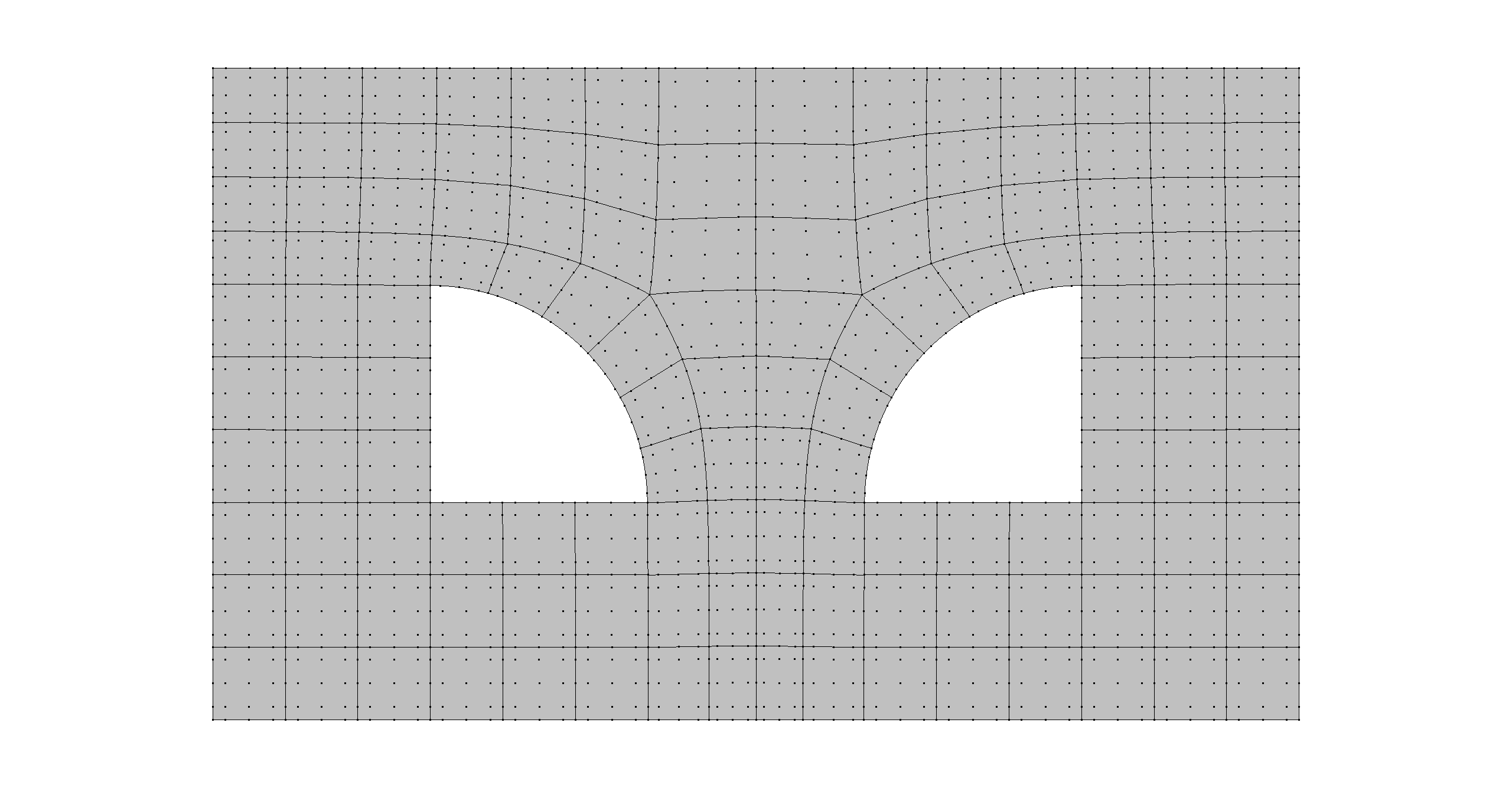}
  \includegraphics[width=\textwidth]{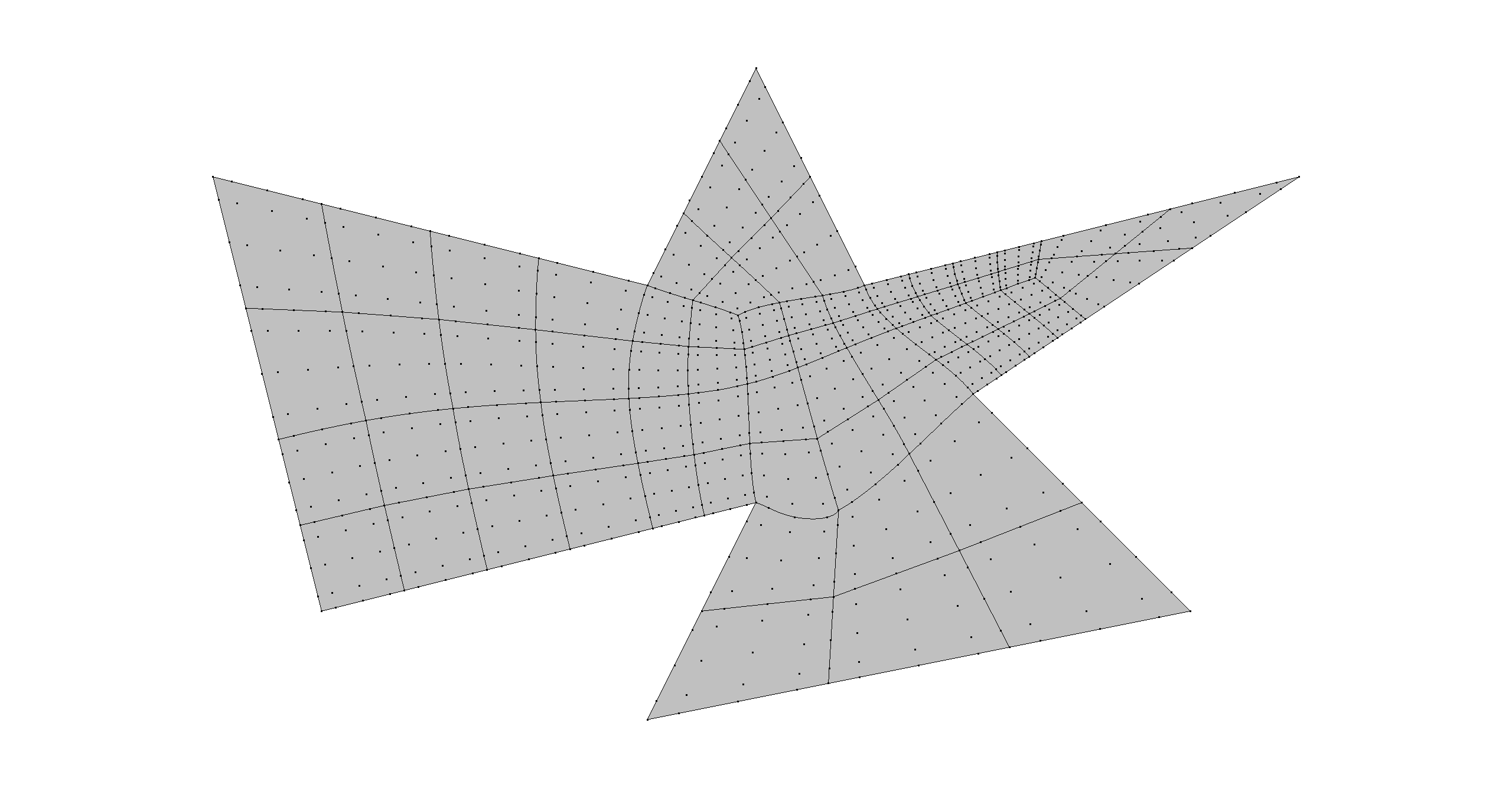}
  \caption{Split meshes obtained on Geometries \rm{I} and \rm{III}.}\label{fig:split-meshes}
\end{figure}

The coarse meshes can be split further according to the user's preference using the isoparametric splitting described in Section~\ref{sec:quad-meshing}.
We show in Fig.~\ref{fig:split-meshes} examples of two split meshes for the quadrilateral decompositions obtained above.

\begin{figure}[htbp]
  \centering
  \includegraphics[width=\textwidth]{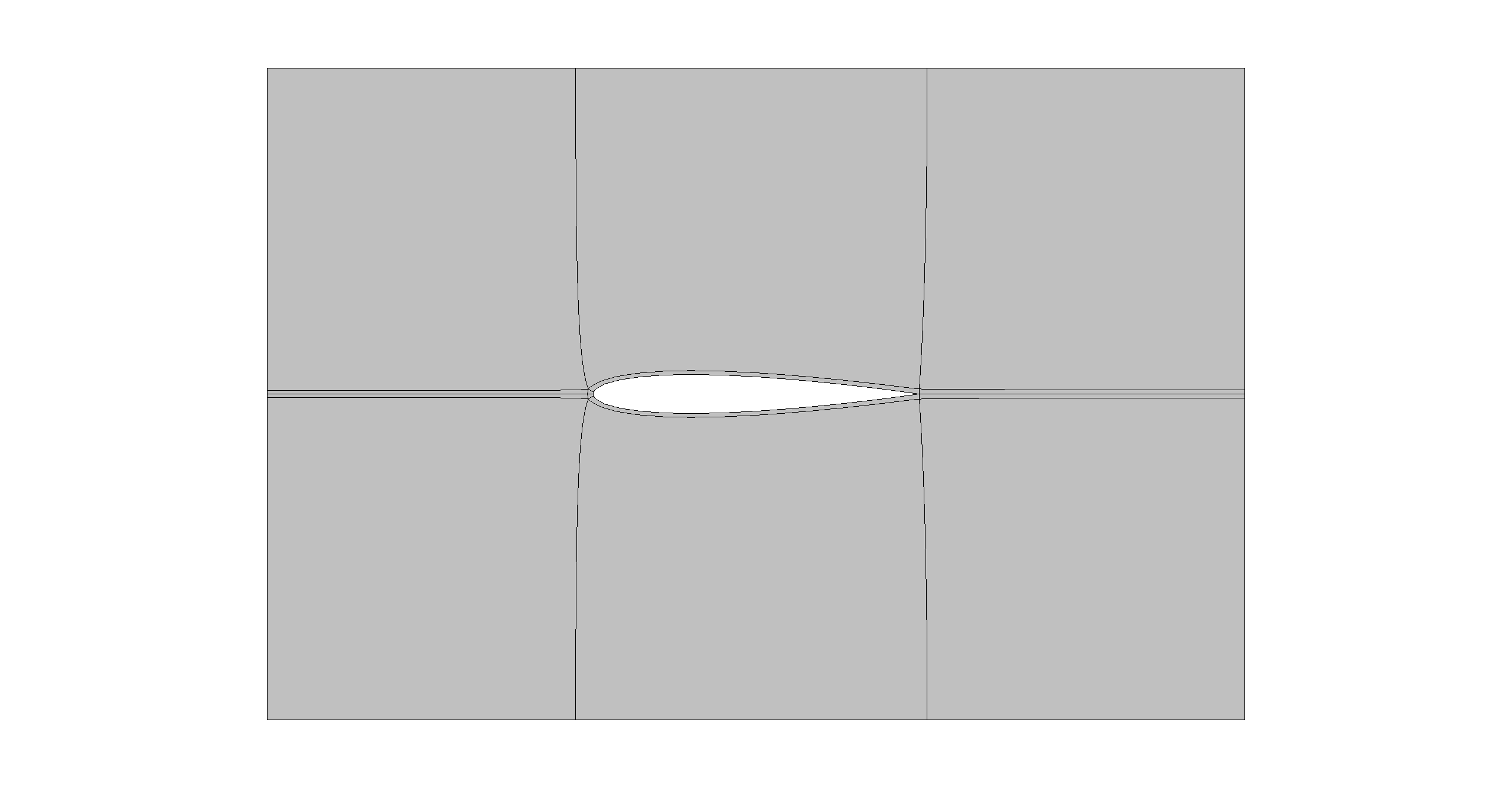}
  \includegraphics[width=\textwidth]{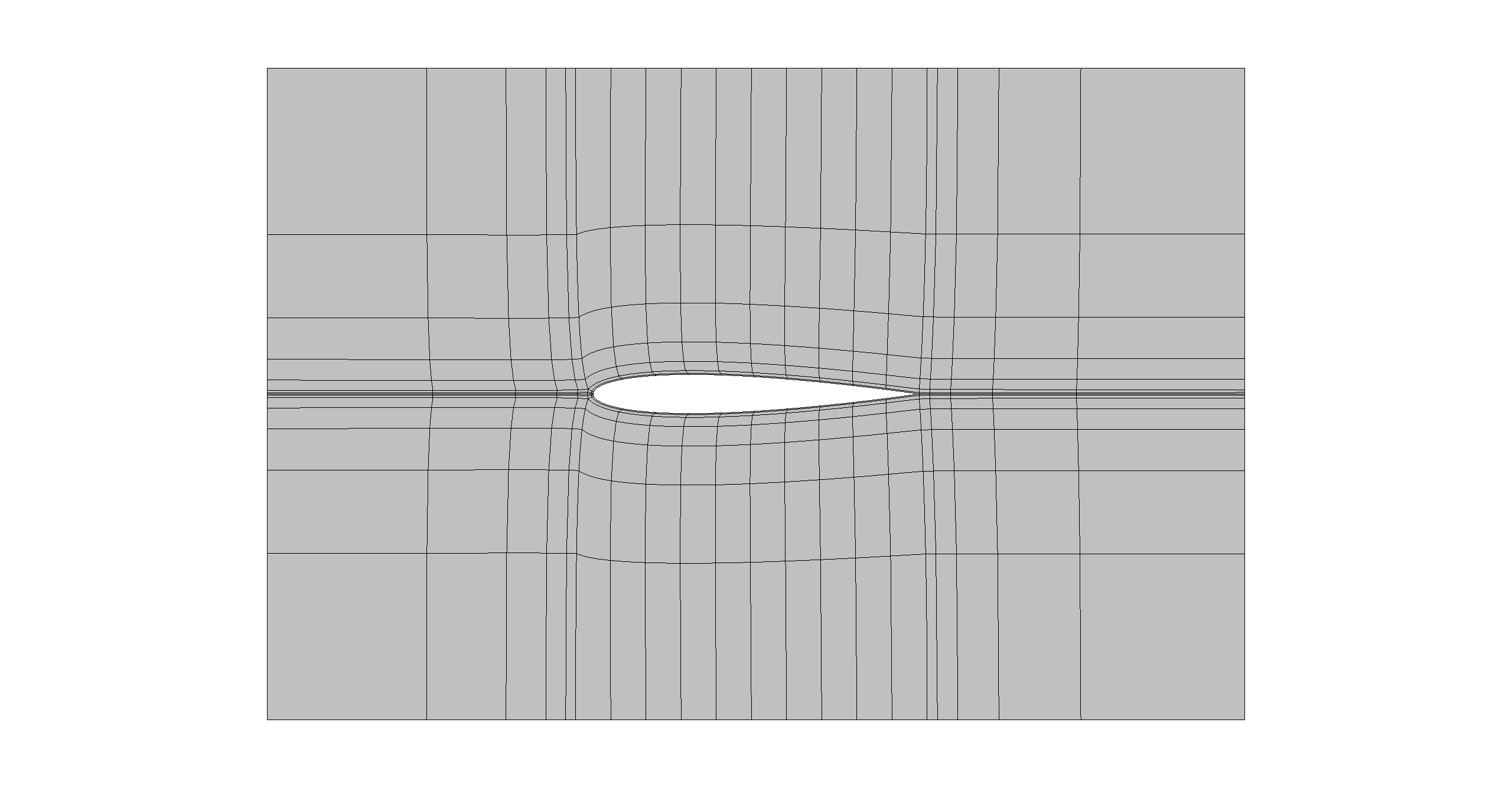}
  \caption{Coarse (top) and split (bottom) quadrilateral meshes obtained on Geometry \rm{IV}, a NACA 0012 profile. High order interior degrees of freedom have been hidden for clearer visualization.}\label{fig:naca-mesh}
\end{figure}

Finally, Example \rm{IV} represents a geometry more relevant for CFD applications, that of a NACA 0012 profile in a rectangular domain.
Since the trailing edge angle is not a multiple of $\pi/2$, a DG discretization was used to compute the guiding field.
With this example, we demonstrate the flexibility of the isoparametric splitting module in \emph{NekMesh} to define different distributions of elements.
This allows the user to easily obtain elements of the required size, e.g.\ high aspect ratio elements in the boundary layer of the airfoil.
Fig.~\ref{fig:naca-mesh} shows the coarse and the split quadrilateral meshes obtained on this NACA 0012 geometry.

\section{Limitations and Extensions}

{\color{black}

The method that we have presented is inspired by cross field methods. Among the possibilities for computing what we are calling the guiding field, we solve a Laplace problem to smoothly propagate the boundary constraints to the interior.
We differ from previous work in two ways:
First, we use high resolution spectral element methods to solve for the guiding field;
second, we do not compute the location of the critical points and separatrices from crosses, but instead from the continuous guiding field.
As a result, the method has advantages over the low order techniques where numerical accuracy is the issue, but retains any disadvantages stemming from the original formulation of the problem.
As a result there are both limitations and areas for further improvement.

The method has the same limitations due to the formulation as do cross field methods.
One such limitation consists in the inability of cross fields to generate singular/critical points (and therefore irregular nodes) in domains where crosses are easily aligned with all nearby boundaries.
Reference~\cite{Fogg2015} documents a modified geometry of the nautilus (of Figs.~\ref{fig:nautilus-streamlines},~\ref{fig:nautilus-solution},~\ref{fig:nautilus-mesh}) with a hole, for which a cross field is unable to generate a valid decomposition.
Fig.~\ref{fig:holed-nautilus-solution} indeed shows that the field is aligned everywhere and that no critical point appears in the domain.
Without critical points, however, a limit cycle is unavoidable and no valid quadrilateral decomposition can be obtained.

\begin{figure}[htbp]
  \centering
  \includegraphics[trim={20cm 0 20cm 0},clip,width=0.49\textwidth]{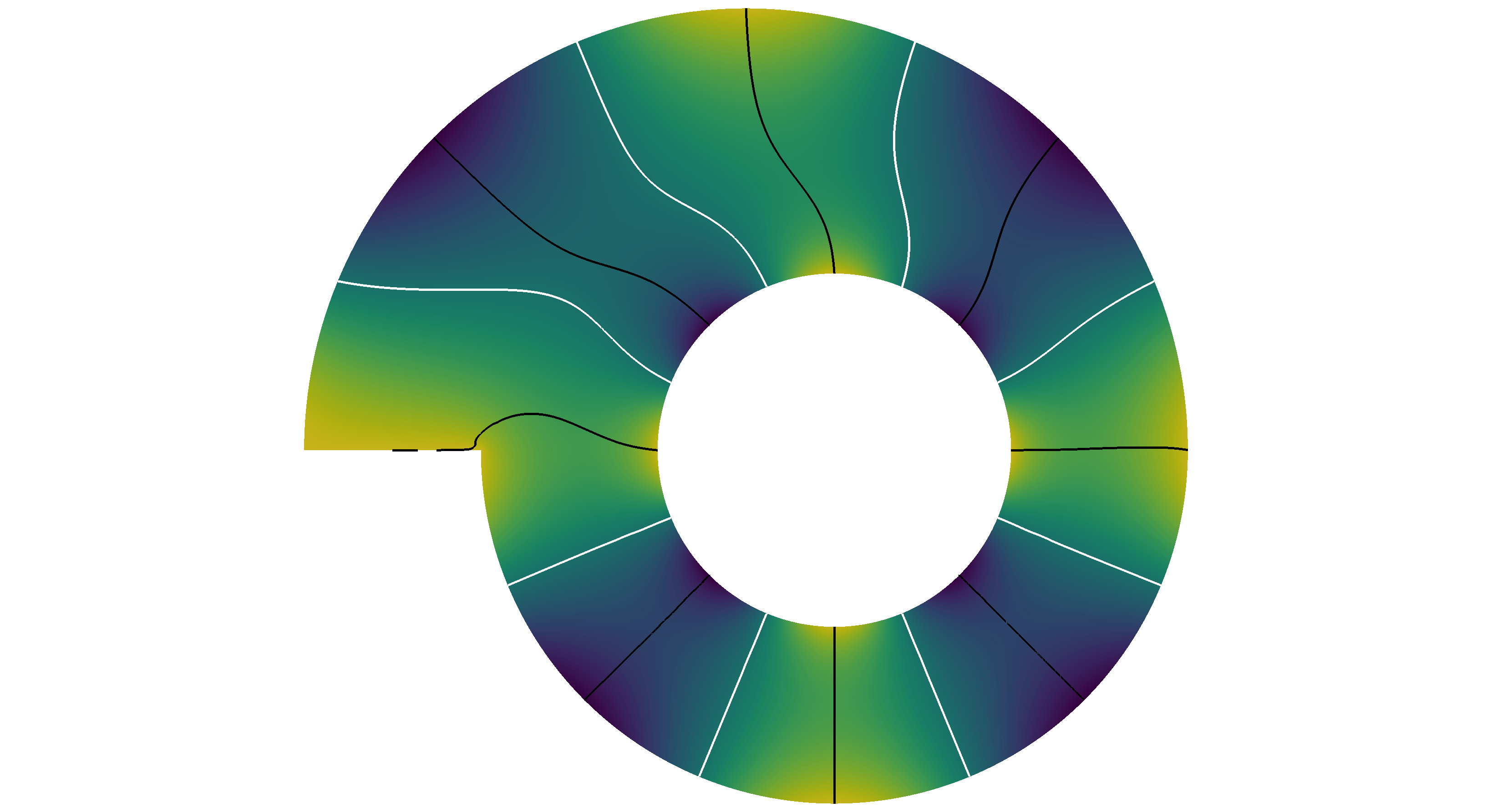}
  \includegraphics[trim={20cm 0 20cm 0},clip,width=0.49\textwidth]{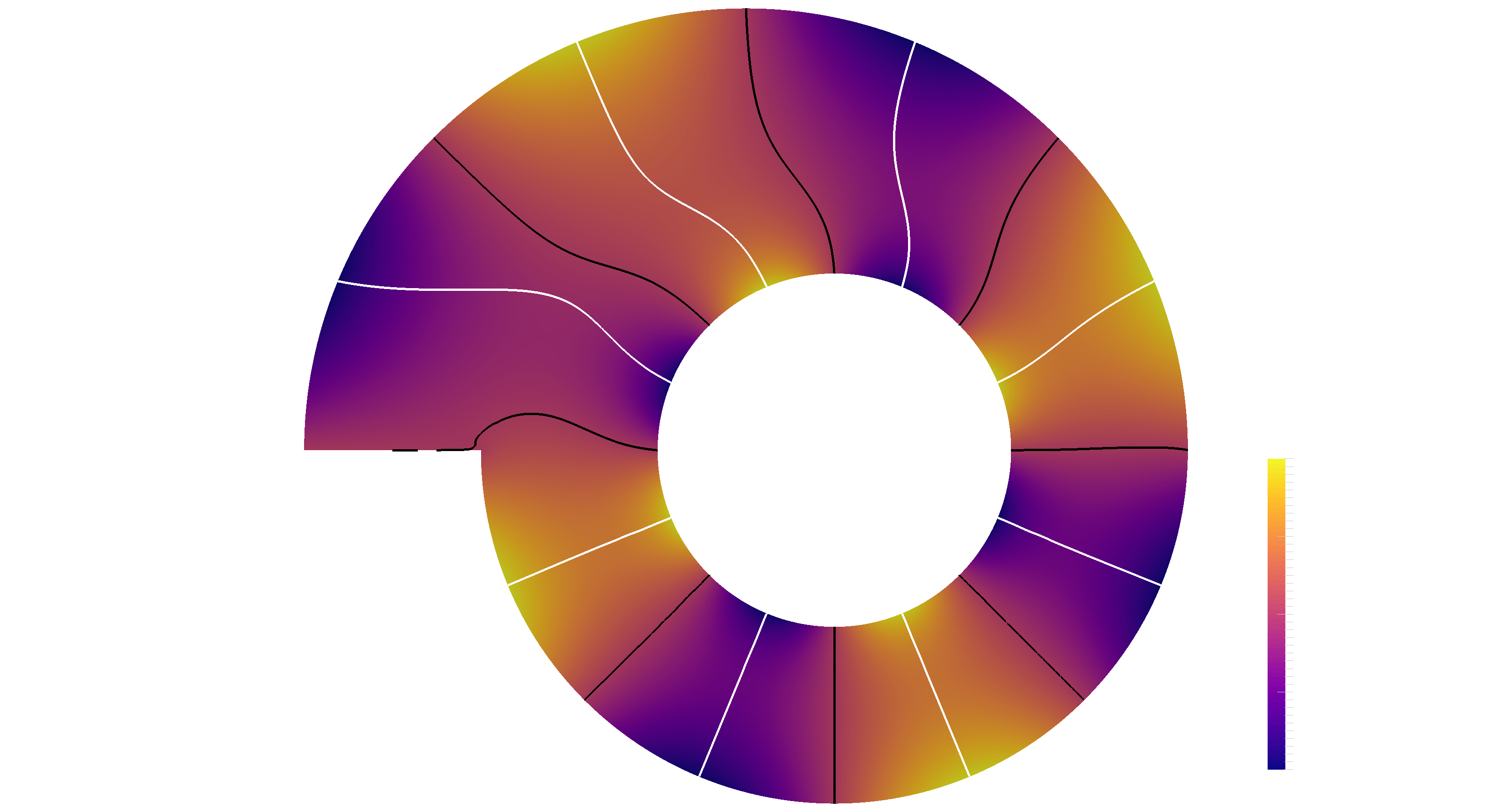}
  \includegraphics[trim={20cm 0 20cm 0},clip,width=0.75\textwidth]{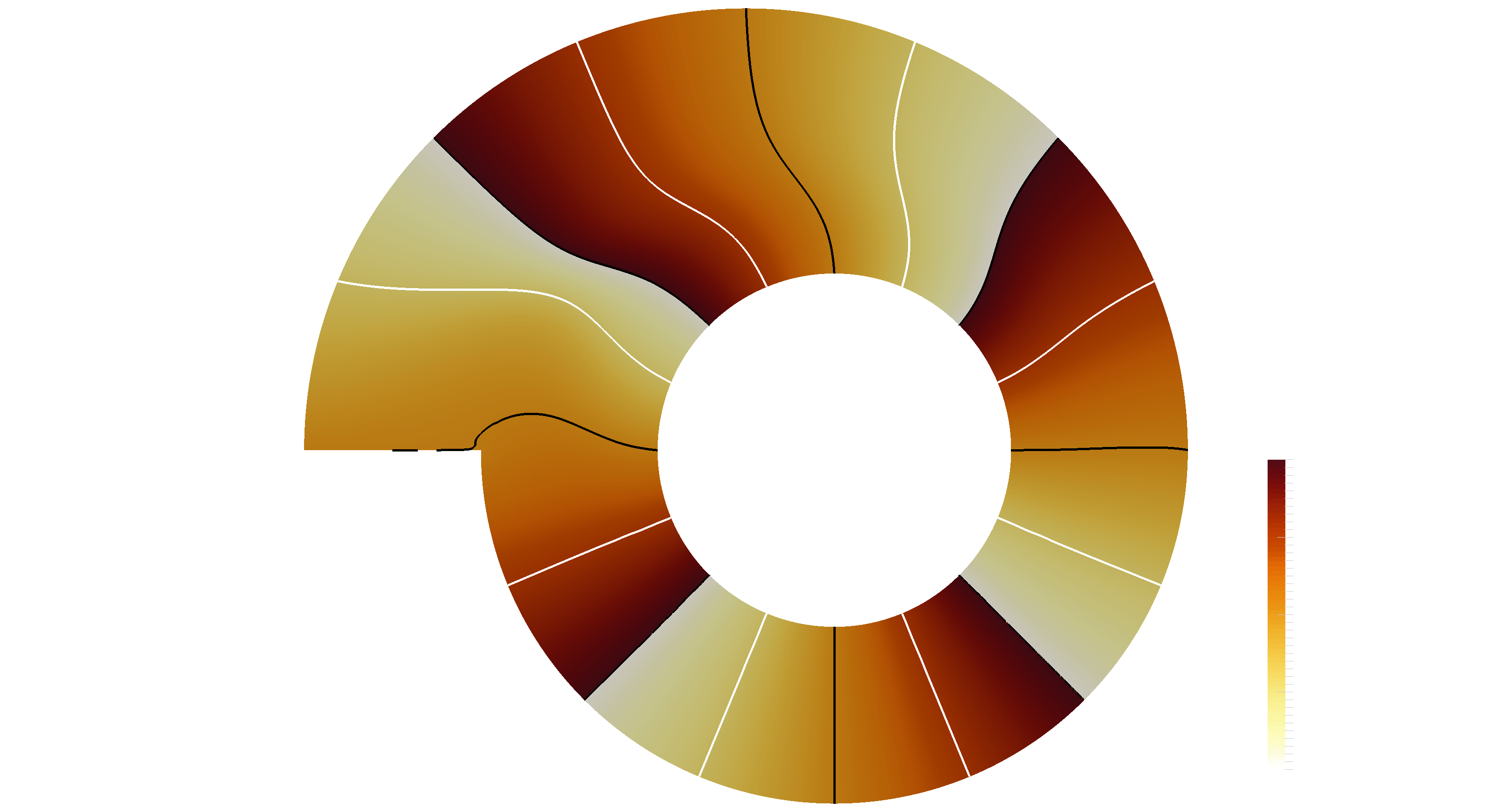}
  \caption{Solution \(u\) (top left) and \(v\) (top right) fields and computed \(\psi \) field (bottom) on the holed nautilus geometry which shows that there are no internal critical points. Isocontours of \(u = 0\) and \(v = 0\) shown in white and black respectively.}\label{fig:holed-nautilus-solution}
\end{figure}

Also, the use of the high order approximations mitigates, but does not eliminate numerical errors.
Therefore, although the nautilus geometry is easily decomposed with the high resolution approximation, it is still possible that spiral streamlines are generated if the accuracy is not high enough. 

The limitations also point the way towards possible improvements to both the formulation and solution of the problem.
Modifications to the guiding field formulation might be used to affect the number and location of critical points and even the existence of a valid quadrilateral decomposition.
Modifications to the algorithm could be made to control errors.

The guiding field could be adjusted by adding a forcing term to \eqref{eq:LaplaceEqns}, thus solving a Poisson rather than a Laplace problem,
which would allow the guiding to be modified field in several ways.
For example, critical points could be added~\cite{2008arXiv0802.2399B,Bunin2006} so that a valid decomposition can be obtained for problems like the holed nautilus. (The degenerate triangle in Fig. \ref{fig:midpoint-division} is also an example of placing an irregular point in the decomposition.)
For the nautilus with hole, a hand sketch shows that adding one 3- and one 5-valent irregular nodes is sufficient.
Existing critical points might also be moved away from boundaries to avoid thin blocks in the decomposition using the idea of control functions used in elliptic grid generation~\cite{Thompsonetal1999}.
Finally, a forcing term would allow one to generate a guiding field that partially aligns with a metric field (incl.\ from a three dimensional surface), such as was demonstrated in \cite{Fogg2015}.

Performance improvements can also be made.
Streamline integration can be performed in reference space~\cite{Coppola2001}, akin to the location and analysis of the critical points.
This would not only be computationally cheaper, it would also allow one to more easily handle geometries with large differences in scales.
Where small scales are present, small elements are typically automatically generated for the background mesh to represent the geometry boundaries accurately.
By performing all of the analysis of the guiding field in reference space, different scales would be automatically handled and better performance might be achieved.
Streamline integration with a symplectic integrator might also be implemented to further mitigate the appearance of numerically created spirals.
Finally, adaptivity and error control could be introduced to further mitigate numerical errors.
Adaptive spectral element methods~\cite{Ekelschot2016} could be used to automatically ensure that the guiding fields errors are acceptably small~\cite{Marcon2019}.
Adaptive ODE time integration methods could further enhance the accuracy of the streamline integration.

An extension of this methodology to volume block decompositions would represent a breakthrough in cross field based techniques.
It would, however, require a reformulation of the problem based on \emph{frame} (three dimensional crosses) arithmetic, which is currently lacking.
Promising work in reference~\cite{Ray2016} points to a nine variable problem to solve.
}

\section{Summary and Conclusions}

We have developed a field guided method to generate quadrilateral meshes for general two dimensional domains.
Inspired by cross field mesh generation methods, the procedure consists of four steps:
The first is to compute the guiding field on an existing triangular mesh.
For that, we use a spectral element method to compute a high resolution approximation to the two Laplace problems for the guiding field.
We use either a continuous or discontinuous Galerkin method depending on the angles formed along the physical boundaries.
If the angles are not multiples of $\pi/2$, the DG formulation can handle the discontinuous boundary conditions in a discretization consistent manner, without the need for \textit{ad hoc} smoothing required by the CG formulation.
{\color{black}Critical} points in the field and their valences are then found by exploiting the local high order polynomial approximations within each element.
Streamlines are computed in the field, again, using the high order solutions rather than crosses defined only at element corners.
Finally, the domain can be cut into elements of the desired size, or used as is.
The procedure is implemented in the open source program \emph{NekMesh}, and uses the companion open source \emph{Nektar++} to compute the guiding fields.

An advantage of the approach is that it can generate meshes with naturally curved quadrilateral elements that do not need to be curved \textit{a posteriori} to eliminate invalid elements.
It avoids low order errors that make it difficult to locate {\color{black}irregular nodes} and accurately trace {\color{black}separatrices}.

\section*{Acknowledgements}

This project has received funding from the European Union's Horizon 2020 research and innovation programme under the Marie Sk\l{}odowska-Curie grant agreement No 675008.
This work was supported by a grant from the Simons Foundation (\#426393, David Kopriva). We thank Ms. Bing Yuan for the plots in Fig. \ref{fig:CrossFieldMeshes}.
The first two authors would also like to thank Prof.~Gustaaf Jacobs of the San Diego State University for his hospitality.

\bibliography{refs}

\end{document}